\documentclass[11pt]{article}
\usepackage[
  margin=2.5cm,
  includefoot,
  footskip=30pt,
]{geometry}

\usepackage[round]{natbib}
\usepackage{amsmath}
\usepackage{symbols}
\usepackage{appendix}
\usepackage{multirow}

\usepackage{xcolor}

\renewcommand{\top}{T}

\long\def\ignore#1{}

\newcommand{\beqn}{\begin{eqnarray}}
\newcommand{\eeqn}{\end{eqnarray}}

\newcommand{\bY}{{\bf Y}}
\newcommand{\bu}{{\bf u}}

\newcommand{\cl}{{\cal L}}
\newcommand{\cb}{{\cal B}}

\newcommand{\sparsity}{\Delta}
\newcommand{\complexity}{r}
\renewcommand{\top}{T}

\newtheorem{assumption}{Assumption}
\newtheorem{subassumption}{Assumption}

\begin{document}

\title{\bf Generalization Error Bounds for Multiclass Sparse Linear Classifiers}
\author{{\bf Tomer Levy} \\ Department of Statistics and Operations Research \\ 
Tel Aviv University
\and 
{\bf Felix Abramovich} \\ Department of Statistics and Operations Research 
\\ Tel Aviv University}

\date{}


\maketitle

\begin{abstract}
We consider high-dimensional multiclass classification by sparse multinomial logistic regression. Unlike binary classification, in the multiclass setup one can think about an entire spectrum of possible notions of sparsity associated with different structural assumptions on the regression coefficients matrix.  We propose a computationally feasible feature selection procedure based on penalized maximum likelihood with convex penalties capturing a specific type of sparsity at hand. In particular, we consider global row-wise sparsity, double row-wise sparsity, and low-rank sparsity, and show that with the properly chosen tuning parameters the derived plug-in classifiers attain the minimax generalization error bounds (in terms of misclassification excess risk) within the corresponding classes of multiclass sparse linear classifiers. The developed approach is  general and can be adapted to other types of sparsity as well. 
\end{abstract}

\section{Introduction}
Classification is a core problem of
statistical and machine learning. One of its main challenges nowadays is high-dimensionality of the data, where the number of features $d$ is of the same order or even larger than the available sample size $n$ (``large $d$, small $n$'' setup) that causes a severe ``curse of dimensionality'' problem. Moreover, the number of classes $L$ may also be large (‘‘large $L$, large $d$, small $n$'' model). A key assumption to handle the ``curse of dimensionality'' is {\em sparsity}. Dimension reduction of the feature domain by selecting a sparse subset of significant features becomes crucial.
\citet{bickel2004fisher} and \citet{fan2008high} showed that even binary classification in high-dimensional setup without a proper feature selection procedure might be as bad as pure guessing.  Feature selection and classification procedures should also be computationally feasible to deal with high-dimensional data.

Although there exists a large amount of statistical and machine learning literature on feature selection in classification, the rigorous theory on the accuracy of resulting classifiers has been mostly developed for the simplest binary case. See  \citet{vapnik2000nature}, \citet{shalev-shwartz2014understanding}, \citet{mohri2018foundations}.

One common strategy for multiclass classification is its reduction to a series of binary classifications. The two most
well-known methods are One-vs-All (OvA), where each class
is compared against all others, and One-vs-One (OvO), where
all pairs of classes are compared to each other. A more direct and appealing strategy is extending binary
classification methods to a  multiclass setup. 
One approach is based on empirical risk minimization (ERM) \citep[e.g.,][]{koltchinskii2002empirical,daniely2012multiclass}. A general crucial drawback of ERM is in minimization of the non-convex 0-1 loss and a common remedy is to replace it by some convex surrogate. \citet{zhang2004statistical}, \citet{chen2006consistency}, \citet{daniely2015multiclass}, \citet{maximov2016tight}, \citet{lei2019datadependent} and \citet{reeve2020optimistic} (see also references therein) investigated the error bounds for various surrogate losses in terms of Rademacher complexity, covering numbers, or Natarajan/graph dimensions. \citet{daniely2015multiclass} compared these results with those for OvA and OvO. However, all the above works do not consider feature selection and to the best of our knowledge, there are no theoretical results for the ERM-based approach
in high-dimensional sparse multiclass setups.

An alternative approach to ERM is plug-in classifiers, where one assumes some model for the underlying unknown probabilities of outcome classes, estimates them from the data and plugs-in estimated probabilities to derive a classification rule. It may be especially useful when one is interested not only in prediction but also in interpretability and inference.
In particular, in this paper we
consider multinomial logistic (linear) classifiers -- one of the mostly used classification tools.
We investigate feature selection in high-dimensional multinomial logistic regression model and the accuracy of the resulting plug-in classifiers under various sparsity scenarios. 

For binary classification the notion of sparsity is naturally associated with the number of
significant features. For linear classifiers it is the number of non-zero entries of a vector of coefficients $\bbeta$. For multiclass case, in contrast, there is a matrix of coefficients $B$ that allows one to consider the entire spectrum of various types of sparsity associated with  different structural assumptions on $B$. \citet{abramovich2021multiclass} studied the most 
straightforward extension of multiclass sparsity measured by the number of non-zero rows of $B$. Such {\em global row-wise}
sparsity corresponds to the assumption that most of features do not affect any class predictions at all. In this paper we present other possible extensions. In particular, we consider {\em double row-wise} sparsity, where it is still assumed that $B$ has a sparse subset of non-zero rows (global sparsity) but, in addition, its non-zero rows are also sparse ({\em local row-wise} sparsity), and the
{\em low-rank} sparsity, where $B$ is assumed to be of a low rank. The latter assumption is associated with the existence of a smaller number of latent variables defining the outcome classes. 

For each considered type of sparsity we propose penalized maximum likelihood feature selection procedures with the corresponding convex penalties and establish the bounds for generalization errors in terms of misclassification excess risk of the resulting multinomial logistic regression classifiers.
The penalties are variations of a celebrated Lasso and its recently developed more general and flexible version Slope 
\citep{bogdan2015slope}. We show that for the proper choice of tuning parameters the derived classifiers attain the optimal (in the minimax sense) generalization errors within the corresponding classes of sparse linear classifiers. The errors can be improved under the additional low-noise condition. 

The proposed approach is  general and can be adapted to other types of sparsity. The machinery for a general form of a sparse multinomial logistic regression classifier is 
developed in Appendix \ref{subsec:appendix_uppper_general}.

The paper is organized as follows.  Section \ref{sec:notations} presents sparse multinomial logistic regression model and some preliminaries.  Section \ref{sec:sparsity} contains the main theoretical results, where we introduce feature selection procedures
for various types of sparsity and derive the error bounds for the resulting misclassification
excess risks. In Section \ref{sec:examples} we illustrate the performance of the developed procedures on a real-data example and compare them with other existing classifiers.
Some concluding remarks are given in Section \ref{sec:conclusion}.
All the proofs are left to the Appendix.

\section{Sparse multinomial logistic regression} \label{sec:notations}
Consider a $d$-dimensional $L$-class classification model:
\be \label{eq:model}
Y|(\bX=\bx) \sim Mult(p_1(\bx),\ldots,p_L(\bx)), 
\;\;\;\;\; \sum_{j=1}^L
p_j(\bx)=1,
\ee
where $\bX \in \mathbb{R}^d$ is a vector of linearly independent features with a marginal probability distribution $\mathbb{P}_X$ on a bounded support $\cX \subseteq \mathbb{R}^d$.  
Let $V=\E(\bX \bX^\top)$ be the second moment matrix of $\bX$.

We consider a multinomial logistic regression model, where it is assumed that
\be \label{eq:multinom}
p_l(\bx)=\frac{\exp(\bbeta_l^T \bx)}{\sum_{k=1}^L\exp(\bbeta_k^T \bx)}, \;\;\; | \bbeta_{l} |_{2} \le R, \;\;\;l=1,\ldots,L,
\ee
Let $B \in \mathbb{R}^{d \times L}$ be the corresponding matrix of regression coefficients with columns $\bbeta_1,\ldots,\bbeta_L$. The model (\ref{eq:multinom}) is not identifiable without an extra constraint on $B$ since shifting each $\bbeta_l$ by the same vector $\bc$ does not affect the probabilities $p_l(\bx)$. 
In this paper we adopt a symmetric zero mean rows constraint $\sum_{l=1}^L \bbeta_l={\bf 0}$
or, equivalently, $B{\bf 1}={\bf 0}_d$. 
Hence, $\bbeta_l$ represents the effects of $\bx$ in the $l$-th class w.r.t. the mean response over all classes on the log-scale:
$$
\bbeta_l^T \bx = \ln \left(\frac{p_l(\bx)}{\sqrt[L]{\prod_{k=1}^L p_k(\bx)}}\right)=\ln p_l(\bx)-\overline{\ln p(\bx)}.
$$
One can evidently choose any other constraint (e.g., $\bbeta_L={\bf 0}$, where the $L$-class is used as the reference one) -- the models will be equivalent but the vectors of coefficients $\bbeta_l$ will have different interpretation.
In particular, the symmetric constraint implies that
the model is invariant to permutations of the classes.

For the considered multinomial logistic regression model (\ref{eq:model})-(\ref{eq:multinom}) the Bayes classifier that minimizes generalized misclassification error (risk) $R(\eta)=P(Y \neq \eta(\bx))$
is a linear classifier 
$\eta^*(\bx)=\arg \max_{1 \leq l \leq L} p_l(\bx)=
\arg \max_{1 \leq l \leq L} \bbeta^T_l \bx$ with the (oracle) misclassification
risk $R(\eta^*)=1-E_{\bX} \max_{1 \leq l \leq L} p_l(\bx)$.

Given a random sample $(\bX_1,Y_1),\ldots,(\bX_n,Y_n)$, we estimate the unknown matrix $B$ from the data and consider the resulting plug-in classifier
${\widehat \eta}_{\widehat B}(\bx)=\arg \max_{1 \leq l \leq L} {\widehat \bbeta}^T_l \bx$. Its conditional misclassification error is $R({\widehat \eta}_{\widehat B})=P\left(Y \neq \widehat{\eta}_{\widehat B}(\bX)|(\bX_1,Y_1),\ldots,(\bX_n,Y_n)\right)$ and its goodness w.r.t. $\eta^*$ is measured by the misclassification excess risk 
$$
\cE(\widehat{\eta}_{\widehat B},\eta^*)=\E R(\widehat{\eta}_{\widehat B})-R(\eta^*).
$$
The goal is to find  $\widehat{B}$ that yields the 
minimal $\cE(\widehat{\eta}_{\widehat B},\eta^*)$.

Consider the log-likelihood function for the multinomial logistic regression model (\ref{eq:model})-(\ref{eq:multinom}):
\be \label{eq:likelihood}
\ell(B)=\sum_{i=1}^n \left\{\bX_i^T B \bxi_i-\ln \sum_{l=1}^L\exp(\bbeta_l^T \bX_i)\right\},
\ee
where $\bxi_i \in \mathbb{R}^L$ is the indicator vector corresponding to $Y_i$ with $\xi_{il}=I\{Y_i=l\}$. One can find the maximal likelihood estimator (MLE) for $B$ by maximizing $\ell(B)$ under the identifiability symmetric constraint $B{\bf 1}={\bf 0}_d$. 
Although the solution is not available in closed form, it can be nevertheless obtained numerically by the fast
iteratively reweighted least squares algorithm \citep{mccullagh1989generalized}.

\ignore{
Assume the following boundedness assumption :
\begin{assumption} \label{as:A}
Assume that there exists $0 < \delta < 1/2$ such that $\delta \leq
p_l(\bx) \leq 1-\delta$
or, equivalently, $|\bbeta_l^T \bx| \leq C^*$ with $C^*=\ln \frac{1-\delta}{\delta}$ for all $\bx \in \cX$ and all $l=1,\ldots,L$. 
\end{assumption}
Assumption \ref{as:A} prevents the conditional variances 
$Var(\xi_l|\bX=\bx)=p_l(\bx)(1-p(\bx))$ to be arbitrarily close to zero, where any MLE-based estimator may fail, and is quite common for GLM  and logistic regression in particular
\citep[e.g.,][]{geer2008highdimensional,alquier2019estimation}.
It also allows one to verify the so-called Bernstein condition which plays a key role in the proofs (see Appendix \ref{subsec:appendix_uppper_general}).
In fact, since Assumption \ref{as:A} is required for all $\bx \in \cX$, without imposing any restrictions on $B$ it essentially implies the boundedness of  $\cX$.
One can replace Assumption \ref{as:A} by other assumptions on $\mathbb{P}_X$. Thus, \citet{alquier2019estimation} showed that the Bernstein condition is also satisfied when ${\bf X}$ is Gaussian.
}

As we have mentioned in the introduction, feature selection is essential for high-dimensional classification.
To perform feature selection we consider a {\em penalized} maximum likelihood estimator of the form:
\be \label{eq:pmle}
\begin{split}
{\widehat B}& =\arg \min_{\widetilde{B}:\widetilde{B}{\bf 1}={\bf 0}_d} \{-\ell(\widetilde{B})+Pen({\widetilde B})\} \\
&=
\arg \min_{\widetilde{B}:\widetilde{B}{\bf 1}={\bf 0}_d} \left\{ \sum_{i=1}^n \left(\ln \sum_{l=1}^L\exp(\bbeta_l^T \bX_i)-\bX_i^T \widetilde{B} \bxi_i\right)+Pen(\widetilde{B}) \right\}
\end{split}
\ee
with a penalty $Pen(\cdot)$ capturing specific sparsity assumptions on $B$. 

\section{Main results}  \label{sec:sparsity}
For binary classification, where a matrix $B$ reduces to a single vector
$\bbeta \in \mathbb{R}^d$, the sparsity is naturally 
characterized by the $l_0$ (quasi)-norm $||\bbeta||_0$ -- the number of non-zero entries of $\bbeta$ \citep[see, e.g.,][]   {abramovich2019highdimensional,chen2021binary}. For the multiclass case there is a wide spectrum of possible ways to extend the notion of sparsity associated with different assumptions on the regression 
coefficients matrix $B$. In this section we consider several of them and derive misclassification excess risk bounds for the resulting multiclass sparse linear classifiers. 

The straightforward approach in (\ref{eq:pmle}) is to use complexity-type penalties that mimic sparsity directly.  However, despite strong theoretical ground \citep[see, e.g.,][]{abramovich2021multiclass}, it is computationally infeasible for high-dimensional data since solving (\ref{eq:pmle}) requires in this case a combinatorial search over all possible models. The goal then is to
find convex surrogates for complexity penalties while preserving their theoretical properties. 

\subsection{Global row-wise sparsity} \label{subsec:global}
The most straightforward extension of notion of sparsity for multiclass classification is {\em global} sparsity, where it is assumed that part of features do not affect any class predictions at all. In terms of the matrix $B$ global sparsity corresponds to the assumption that $B$ has a ``small'' number of non-zero rows (global row-wise sparsity). Such type of sparsity was studied in \citet{abramovich2021multiclass}
and in this subsection we review their main results (generalizing for a anisotropic $\bX$) in order to extend them afterwards to other, finer types of sparsity. 

Let $r_B$ be the number of non-zero rows of $B$.
To capture the global sparsity \citet{abramovich2021multiclass} proposed  to use a complexity penalty on the number of non-zero rows of $B$ in (\ref{eq:pmle}).

Let \(\cM = \{ B \in \R^{d \times L} : B{\bf 1}={\bf 0} \}\) be the set of regression matrices satisfying the symmetric constraint,
$\cM(d_0)=\{B \in \cM: r_B \leq d_0\}$ be its subset of  $d_0$-globally row-wise sparse matrices
and $\cC_L(d_0)=\{\eta(\bx)=\arg \max_{1 \leq l \leq L} \bbeta_l^T \bx: B \in \cM(d_0)\}$
be the set of $d_0$-sparse linear $L$-class classifiers. Define the penalized maximum likelihood estimator $\hat{B}$ 
\be \label{eq:penalized}
\widehat{B}=\arg \min_{\widetilde{B} \in \cM} \left\{-\ell(\widetilde{B})+Pen(r_{\widetilde B})\right\}
\ee
with the complexity penalty of the form
\be \label{eq:optpen}
Pen(r_{\widetilde B})=C_1~r_{\widetilde B}(L-1)+C_2~ r_{\widetilde B}\ln\left(\frac{d e}{r_{\widetilde B}}\right)
\ee 
for some positive constants $C_1$ and $C_2$.

\citet{abramovich2021multiclass} showed that for the bounded \(\cX\), 
\be \label{eq:upper}
\sup_{\eta^* \in \cC_L(d_0)}
\cE(\widehat{\eta}_{\widehat B},\eta^*) \leq C\sqrt{\frac{d_0 (L-1)+   
d_0 \ln \left(\frac{de}{d_0}\right)}{n}}
\ee
for some $C>0$ simultaneously for all $1 \leq d_0 \leq \min(d,n)$,
and that the bound in (\ref{eq:upper}), up to a probably different constant, is also the minimax for $\cC_L(d_0)$.

Misclassification excess risk bounds (\ref{eq:upper}) show that there is a phase transition between small and large number of classes. For $L \leq 2+\ln(d/d_0)$, the multiclass effect is not yet manifested and the minimax
misclassification excess risk over the set of $d_0$-sparse linear
classifiers is of the order
$\sqrt{\frac{d_0}{n}\ln\left(\frac{de}{d_0}\right)}$ regardless of $L$. Note that $d_0\ln\left(\frac{de}{d_0}\right) \sim \ln\binom{d}{d_0}$ which is the log of the number of all possible models of size $d_0$. For larger $L$, the  risk is of the order $\sqrt{\frac{d_0(L-1)}{n}}$, where $d_0(L-1)$ is the overall number of estimated parameters in the given model of size $d_0$, and does not depend on $d$.
For $L>n/d_0$
the number of parameters in the (true) model becomes
larger than the sample size and consistent classification is
evidently impossible.

\ignore{
\begin{theorem}[Theorem 2 of \citet{abramovich2021multiclass}] \label{th:lower}
Consider a $d_0$-sparse agnostic multinomial logistic regression
model (\ref{eq:model})-(\ref{eq:multinom}), where
$2 \leq  d_0 \ln(\frac{de}{d_0}) \leq n$ and $d_0(L-1) \leq n$. Then,
$$
\inf_{\widetilde \eta}\sup_{\eta^* \in \cC_L(d_0),~\mathbb{P}_X}
\cE(\widetilde{\eta},\eta^*) \geq C_2 {\sqrt \frac{d_0 (L-1)+d_0
\ln \left(\frac{de}{d_0}\right)}{n}}
$$
for some $C_2>0$, where the infimum is taken over all classifiers $\tilde{\eta}$ based on the data $(\bX_i,Y_i),\;i=1,\ldots,n$.
\end{theorem}
}
Classification is mostly challenging at points,
where it is difficult to distinguish the most likely class from
others, that is, at those $\bx \in \cX$, where the largest probability $p_{(1)}(\bx)$ is close to the second largest $p_{(2)}(\bx)$.
The misclassification error bounds (\ref{eq:upper}) may be then improved under the additional multiclass extension of the low-noise (aka Tsybakov) condition \citep{tsybakov2004optimal}:
\begin{assumption} \label{as:B}
Consider the multinomial logistic regression model (\ref{eq:model})-(\ref{eq:multinom}) and assume that there exist
$C>0, \alpha \geq 0$ and $h^*>0$ such that for all $0<h \leq h^*$,
$$
P\left(p_{(1)}(\bX)-p_{(2)}(\bX) \leq h\right) \leq C h^{\alpha}.
$$
\end{assumption}
\noindent
Assumption \ref{as:B} implies that with high probability (depending on the parameter $\alpha$) the most likely  class is sufficiently distinguished from others. The two extreme cases are $\alpha=0$ and
$\alpha=\infty$. The former does not impose any assumption on the noise, while the latter assumes the existence of a hard margin of size $h^*$ separating $p_{(1)}(\bx)$ and $p_{(2)}(\bx)$.

\citet{abramovich2021multiclass} proved that under the additional low-noise Assumption \ref{as:B} the misclassification excess risk bound (\ref{eq:upper}) of $\widehat{\eta}_{\widehat B}$ can be indeed improved: 
\be \label{eq:upper_noise}
\sup_{\eta^* \in \cC_L(d_0)}\cE(\widehat{\eta}_{\widehat B},\eta^*) \leq \left(C~ \frac{d_0 (L-1)+   
d_0 \ln \left(\frac{de}{d_0}\right)}{n}\right)^{\frac{\alpha+1}{\alpha+2}}
\ee
for all $1 \leq d_0 \leq \min(d,n)$ and all $\alpha \geq 0$.
Note that the proposed classifier $\widehat{\eta}_{\widehat B}$ is inherently adaptive to both sparsity
$d_0$ and noise level $\alpha$.

As we have mentioned above, 
solving for $\widehat{B}$ in (\ref{eq:penalized}) requires a combinatorial search over all possible $2^d$ models that
makes it computationally infeasible for large $d$. One should apply convex relaxation techniques to replace the original complexity penalty (\ref{eq:optpen}) by some convex surrogate.

The well-known examples of convex surrogates are 
celebrated Lasso, where the $l_0$-norm in the complexity penalty is replaced by the $l_1$-norm
norm, and its recently developed more general variation Slope that uses a {\em sorted} $l_1$-type norm \citep{bogdan2015slope}. 
Lasso and Slope estimators have been intensively studied in
the last decade in various regression setups \citep[see e.g.,][]{geer2008highdimensional,bickel2009simultaneous,su2016slope,bellec2018slope,abramovich2019highdimensional,alquier2019estimation}.
\citet{abramovich2019highdimensional} and \citet{abramovich2021multiclass} applied logistic Lasso and Slope classifiers for classification.

To capture a global row-wise sparsity for multinomial logistic regression, \citet{abramovich2021multiclass} considered a {\em group} version of multinomial logistic Slope defined as follows.
Let
\be \label{eq:groupSlope}
\widehat{B}_{gS} =\arg \min_{\widetilde B \in \cM}
\left\{\frac{1}{n}\sum_{i=1}^n\left(\ln\left(\sum_{l=1}^L \exp(\widetilde{\bbeta}_l^T \bX_i)\right)-\bX_i^T\widetilde{B}\bxi_i\right)+ \sum_{j=1}^d \lambda_j |{\widetilde B}|_{(j)} \right\},
\ee
where $|\widetilde{B}|_{(1)} \geq \ldots \geq
|\widetilde{B}|_{(d)}$ are the descendingly ordered $l_2$-norms of rows of $\widetilde{B}$
and $\lambda_1 \geq \ldots \geq \lambda_d > 0$ are tuning parameters, and define $\widehat{\eta}_{gS}(\bx)=\arg\max_{1 \leq l \leq L}\widehat{\bbeta}_{gS,l}^T \bx$. 
Multinomial logistic group Lasso classifier $\widehat{\eta}_{gL}$ is a particular case of $\widehat{\eta}_{gS}$ corresponding to equal $\lambda_j$'s in  (\ref{eq:groupSlope}).

The identifiability symmetric constraint ${\widetilde B} \in \cM$ is, in fact, unnecessary in (\ref{eq:groupSlope}) since 
unlike the complexity penalty in (\ref{eq:penalized}), the solution of (\ref{eq:groupSlope}) is identifiable without any additional constraint.
Moreover, since the unconstrained log-likelihood (\ref{eq:likelihood}) satisfies $\ell(\widetilde{\bbeta}_1,\ldots,\widetilde{\bbeta}_L)=\ell(\widetilde{\bbeta}_1-{\bf c},\ldots,\widetilde{\bbeta}_L-{\bf c})$ for any  
vector ${\bf c} \in \mathbb{R}^d$, it can be always improved by taking $\hat{c}_j=\arg\min_{c_j} \sum_{l=1}^L (\widetilde{B}_{jl}-c_j)^2$, that is, for 
$\hat{c}_j=\bar{B}_{j\cdot}$. 
Hence, the unconstrained solution of (\ref{eq:groupSlope})  
will inherently have zero mean rows.

As usual for convex relaxation, one needs some (mild) extra conditions on the design. Assume that all $X_j$ are scaled, i.e. $EX_j^2=1,\;j=1,\ldots,d$. For a given matrix $A \in \cM$ let
\(\Pi_{d_0}(A)\) be its $d_0$-sparse projection, i.e. the matrix with at most $d_0$ nonzero rows closest to $A$ in the Frobenius norm.  

\ignore{
\begin{subassumption} \label{as:Dgs} Let $\nu_{gS}(d_0)$ be the minimal $d_0$-sparse eigenvalue of $V$, i.e.
$$
\nu_{gS}(d_0) =\inf_{A \in \cM(d_0): A \neq 0_{d \times L}} \frac{||V^\frac{1}{2}A||^2_F}{||A||^2_F}
$$
and assume that $\nu_{gS}(d_0) >0$.
\end{subassumption}
}

\begin{subassumption} \label{as:Dgs}
    Assume that 
    \[
    \nu_{gS}(d_0)=\inf_{A \in \cM: A \neq 0_{d \times L}} \frac{\|V^{\frac{1}{2}} A\|^2_{F}}{\| \Pi_{d_0}(A) \|^2_F} > 0,
    \]
\end{subassumption}
In fact, one immediately realizes that \(\Pi_{d_0}(A)\) keeps $d_0$ rows of $A$ with the largest $l_2$-norms and zeroes other rows. Hence,
\(\| \Pi_{d_0}(A) \|^2_F=\sum_{j=1}^{d_0} |A|^2_{(j)}\).

Such or similar types of conditions are common for convex relaxation methods \citep[see ][Section 8 for discussion]{bellec2018slope}.

Let $||A||_{gS}=\sum_{j=1}^d \lambda_j |A|_{(j)}$ be the group Slope norm of a matrix $A$.
The following theorem provides an upper bound for misclassification excess risk of the group Slope classifier extending the results of \citet{abramovich2021multiclass} to anisotropic design. In addition, it provides also the upper bounds for the integrated prediction error $\sum_{l=1}^L \E\|(\widehat{\bbeta}_{gS,l}-\bbeta_l)^T\bx\|^2_{L_2(\mathbb{P}_X)}=E\|V^{\frac{1}{2}}(\widehat{B}_{gS}-B)\|^2_F$ and the estimation error of the regression coefficients matrix $B$ w.r.t. the group Slope norm $\E ||\widehat{B}_{gS}-B||_{gS}$ :
\begin{theorem} \label{th:gs_low_noise_classifier} Consider a $d_0$-globally row-sparse multinomial logistic regression (\ref{eq:model})-(\ref{eq:multinom}), where $\bX_j$'s are scaled to have $EX_j^2=1,\;j=1,\ldots,d$. Apply the multinomial logistic sparse group Slope classifier (\ref{eq:sgS}) with $\lambda_j$'s satisfying 
\be \label{eq:weightsgs}
\max_{1 \leq j \leq d} \frac{\sqrt{L + \ln(d / j)}}{\lambda_{j}} \le C_0 \sqrt{n},
\ee
where the constant \(C_0\) is derived from \citet{abramovich2021multiclass}. Then, under Assumptions \ref{as:B}-\ref{as:Dgs},
$$
\sup_{\eta^* \in \cC_L(d_0)}\cE(\widehat{\eta}_{gS},\eta^*) 
\le \left(\frac{C}{\nu_{gS}(d_0)} \sum_{j=1}^{d_0} 
\frac{\lambda_j}{\sqrt{j}}\right)^{\frac{2(\alpha+1)}{\alpha+2}}.
$$
In addition, 
$$
\sup_{B \in \cM(d_0)} \E\|V^{\frac{1}{2}}(\widehat{B}_{gS}-B)\|^2_F \leq \frac{C_1}{\nu_{gS}(d_0)} \left( \sum_{j=1}^{d_0} 
\frac{\lambda_j}{\sqrt{j}} \right)^2
$$
and
$$
\sup_{B \in \cM(d_0)} \E ||\widehat{B}_{gS}-B||_{gS} \leq \frac{C_2}{\nu_{gS}(d_0)} \left( \sum_{j=1}^{d_0} 
\frac{\lambda_j}{\sqrt{j}} \right)^2
$$
\end{theorem}
In particular, setting
$$
\lambda_j = \frac{1}{C_0} \sqrt{\frac{L+\ln(d/j)}{n}},
$$ 
the misclassification excess risk of the multinomial logistic group Slope classifier $\widehat{\eta}_{gS}$ is of the minimax order (\ref{eq:upper_noise}): 
\begin{corollary} \label{cor:groupslope}
Apply Theorem \ref{th:gs_low_noise_classifier} with
$$
\lambda_j = \frac{1}{C_0} \sqrt{\frac{L+\ln(d/j)}{n}}.
$$
Then,
under Assumptions \ref{as:B}-\ref{as:Dgs},
\be \label{eq:gSerror}
\sup_{\eta^* \in \cC_L(d_0)}\cE(\widehat{\eta}_{gS},\eta^*) 
\le \left(\frac{C}{\nu_{gS}(d_0)}~\frac{d_0(L-1)+d_0 \ln\left(\frac{de}{d_0}\right)}{ n}\right)^{\frac{\alpha+1}{\alpha+2}}.
\ee

Furthermore,
$$
\sup_{B \in \cM(d_0)} \E\|V^{\frac{1}{2}}(\widehat{B}_{gS}-B)\|^2_F \leq \frac{C_1}{\nu_{gS}(d_0)}~\frac{d_0(L-1)+d_0 \ln\left(\frac{de}{d_0}\right)}{n}
$$
and
$$
\sup_{B \in \cM(d_0)} \E ||\widehat{B}_{gS}-B||_{gS} \le \frac{C_2}{\nu_{gS}(d_0)}~\frac{d_0(L-1)+d_0 \ln\left(\frac{de}{d_0}\right)}{n}
$$
\end{corollary} 
Note that $\widehat{\eta}_{gS}$ is inherently adaptive to $d_0$ and $\alpha$.

Similarly, the multinomial logistic group Lasso classifier
$\widehat{\eta}_{gL}$ with a (constant) $\lambda = \frac{1}{C_0} \sqrt{\frac{L +\ln d}{n}}$ is sub-optimal (up to the log-factor):
$$
\sup_{\eta^* \in \cC_L(d_0)}\cE(\widehat{\eta}_{gL},\eta^*) 
\le \left(\frac{C}{\nu_{gS}(d_0)}~ \frac{d_0(L-1)+d_0 \ln d}{n}\right)^{\frac{\alpha+1}{\alpha+2}}.
$$

We consider now other possible types of sparsity for multiclass case and derive the corresponding generalization error bounds.

\subsection{Double row-wise sparsity} \label{subsec:double}
We show that the misclassification excess risks bounds for a global row-wise sparsity can be improved under a finer row-wise sparsity structure. 
Namely, assume that even each significant feature is involved in only part of probabilities $p_l$'s. In terms of the matrix $B$
it implies the additional sparsity assumption on its non-zero rows in the usual $l_0$-norm sense, i.e., local row-wise sparsity.

For a given matrix $B$, let \(\cJ(B) = \{j_{1}, \ldots, j_{r_B}\}\) be the set of indices of its non-zero rows. Consider a set of {\em double} (global and local) row-wise sparse matrices 
$\cM(d_0,{\bf m})=\{B \in \cM: |\cJ(B)| \leq d_0;\; ||B_{j\cdot}||_0 \leq m_{j}, \;j \in \cJ(B)\}$ and 
the corresponding set of $(d_0,{\bf m})$-sparse linear $L$-class classifiers
$\cC_L(d_0,{\bf m})=\{\eta(\bx)=\arg \max_{1 \leq l \leq L} \bbeta^T_l\bx: B \in \cM(d_0,{\bf m})\}$.

To capture a double row-wise sparsity one should 
impose penalties on both the number of non-zero rows of $B$ and 
on the numbers of their non-zero entries. A natural
convex surrogate in this case is 
a multinomial logistic {\em sparse} group Slope estimator of $B$ defined as follows:
\begin{equation}\label{eq:sgS}
\begin{split}
\widehat{B}_{sgS} =\arg \min_{\widetilde B \in \cM}
& \left\{\frac{1}{n}\sum_{i=1}^n\left(\ln\left(\sum_{l=1}^L \exp(\widetilde{\bbeta}_l^T \bX_i)\right)-\bX_i^T\widetilde{B}\bxi_i\right) \right. \\
& \left. + \sum_{j=1}^d \lambda_j |{\widetilde B}|_{(j)}+
\sum_{j=1}^d\sum_{l=1}^L \kappa_l |\widetilde{B}|_{j(l)}
 \right\},
\end{split}
\ee
where $|\widetilde{B}|_{(1)} \geq \ldots \geq |\widetilde{B}|_{(d)}$ are the descendingly ordered $l_2$-norms of the rows of $\widetilde{B}$, $|\widetilde{B}|_{j(1)} \geq \ldots \geq|\widetilde{B}|_{j(L)}$ are the descendingly ordered absolute values of entries of its $j$-th row, and 
$\lambda_1 \geq \ldots \geq \lambda_d > 0$ and $\kappa_1 \geq \ldots \geq \kappa_L > 0$ are tuning parameters.
The additional last term in the penalty in (\ref{eq:sgS}) yields sparsity of non-zero rows. Sparse group Slope essentially combines group Slope on the row's norms with
usual Slope within each row.

The multinomial logistic sparse group Slope classifier  is
$\widehat{\eta}_{sgS}(\bx)=\arg \max_{1 \leq l \leq L}
\widehat{\bbeta}_{sgS,l}^T \bx$.
Multinomial logistic sparse group Lasso  classifier $\widehat{\eta}_{sgL}$ \citep[see][]{friedman2010regularization,vincent2014sparse} is its particular case  
with identical $\lambda_j$'s and $\kappa_l$'s in (\ref{eq:sgS}).

\ignore{
\begin{subassumption} \label{as:Dsgs} Let $\nu_{sgS}(d_0, \bm)$ be the minimal $(d_0,\bm)$-sparse eigenvalue of $V$, i.e.
$$
\nu_{sgS}(d_0,\bm) =\inf_{A \in \cM(d_0,\bm):A \neq 0_{d \times L}} \frac{||V^\frac{1}{2}A||^2_F}{||A||^2_F}
$$
and assume that $\nu(d_0,\bm) >0$.
\end{subassumption}
}

\ignore{
\begin{subassumption} \label{as:Dsgs}
Assume that for all matrices \(A \in \R^{d \times L}\) such that
$||A||_{sgS} = 10\sum_{j=1}^{d_0} \frac{\lambda_{j}}{\sqrt{j}} + 10\sum_{j=1}^{d} \sum_{l=1}^{m_j} \frac{\kappa_{l}}{\sqrt{l}}$ and
$||V^{\frac{1}{2}}A||^2_F \leq 1$,
    \[
    \nu_{sgS}(d_0, \bm) = \inf_A \frac{||V^{\frac{1}{2}} A||^2_F}{||\Pi_{d_0}(A)||^2_F} > 0.
    \]
    \end{subassumption}
}

Let $||A||_{sgS}=\sum_{j=1}^{d} \lambda_{j} |A|_{(j)} + \sum_{j=1}^{d} \sum_{l=1}^{L} \kappa_{l} |A|_{j(l)}$ be the sparse group Slope norm of a matrix $A \in \mathbb{R}^{d \times L}$.
The following theorem provides an upper bound for  misclassification excess risk of $\widehat{\eta}_{sg S}$:
\begin{theorem} \label{th:sgs_low_noise_classifier} Consider a $(d_0,{\bf m})$-sparse multinomial logistic regression (\ref{eq:model})-(\ref{eq:multinom}) with scaled $X_j$'s. Apply the multinomial logistic sparse group  classifier (\ref{eq:sgS}) with 
$\lambda_j$'s and $\kappa_l$'s satisfying \(\kappa_{L} \ge \sqrt{\frac{\pi}{2}}\frac{2880}{7} \frac{1}{\sqrt{n}}\) and
\be \label{eq:weightssgs}
\max_{1 \leq j \leq d} \frac{\sqrt{2 \sum_{l=1}^{L}\frac{1}{l}\left(\frac{Le}{l}\right)^{l}e^{-C^2 n l \kappa_{l}^{2}}+2\log\left(\frac{de}{j}\right)}}{\lambda_{j}} \le \frac{7}{1440 C_0} \sqrt{\frac{2}{\pi}}\sqrt{n},
\ee
\(C = \sqrt{\frac{2}{\pi}}\frac{7}{2880}\) and \(C_0\) is derived in the proof.
Then, under Assumptions \ref{as:B}-\ref{as:Dgs},
$$
\sup_{\eta^* \in \cC_L(d_0, {\bf m})}\cE(\widehat{\eta}_{sgS},\eta^*) 
\le \left(\frac{C}{\nu_{gS}(d_0)} \left(\sum_{j=1}^{d_0} 
\frac{\lambda_j}{\sqrt{j}} + \sqrt{\sum_{j=1}^{d_0} 
\left(\sum_{l=1}^{m_j} \frac{\kappa_l}{\sqrt{l}}\right)^2} \right) \right)^{\frac{2(\alpha+1)}{\alpha+2}}.
$$
In addition,
$$
\sup_{B \in \cM(d_0)} \E\|V^{\frac{1}{2}}(\widehat{B}_{sgS}-B)\|^2_F \leq \frac{C_1~}{\nu_{gS}(d_0)} \left(\sum_{j=1}^{d_0} 
\frac{\lambda_j}{\sqrt{j}} + \sqrt{\sum_{j=1}^{d_0} 
\left(\sum_{l=1}^{m_j} \frac{\kappa_l}{\sqrt{l}}\right)^2} \right)^2
$$
and
$$
\sup_{B \in \cM(d_0,{\bf m})}  \E||\widehat{B}_{sgS}-B||_{sgS} \leq 
\frac{C_2~}{\nu_{gS}(d_0)} \left(\sum_{j=1}^{d_0} 
\frac{\lambda_j}{\sqrt{j}} + \sqrt{\sum_{j=1}^{d_0} 
\left(\sum_{l=1}^{m_j} \frac{\kappa_l}{\sqrt{l}}\right)^2} \right)^2.
$$
\end{theorem}         
The proof is given in the Appendix \ref{sec:appendix_upper}.

In particular, take $\lambda_j=c_1 \sqrt{\frac{\ln\left(de/j\right)}{n}},\;j=1,\ldots,d$ and
$\kappa_l=c_2\sqrt{\frac{\ln\left(Le/l\right)}{n}},\;l=1,\ldots,L$ with $c_1=\frac{1440 C_0 \sqrt{\pi}}{7}$ and $c_2=\frac{2880 \sqrt{\pi}}{7}$. By a straightforward calculus one can verify that these $\lambda_j$'s and $\kappa_l$'s satisfy the condition 
(\ref{eq:weightssgs}), and Theorem \ref{th:sgs_low_noise_classifier} then implies:
\begin{corollary} \label{cor:complexity_bound}
Apply Theorem \ref{th:sgs_low_noise_classifier} with
\be \label{eq:sgslope_parameters}
\lambda_j=c_1\sqrt{\frac{\ln\left(de/j\right)}{n}} \;\;\;{\rm and}\;\;\;
\kappa_l=c_2 \sqrt{\frac{\ln\left(Le/l\right)}{n}},
\ee
where $c_1=\frac{1440 \sqrt{2\pi} C_0}{7},\; c_2 =  \frac{2880\sqrt{\pi}}{7}$ and $C_0$ is given in Lemma \ref{lem:rad_truncated}.
Then,
under Assumptions \ref{as:B}-\ref{as:Dgs}, 
\be \label{eq:sgSerror}
\sup_{\eta^* \in \cC_L(d_0, {\bf m})}\cE(\widehat{\eta}_{sgS},\eta^*) 
\le \left(\frac{C}{\nu_{gS}(d_0)}~\frac{d_0 \ln\left(\frac{de}{d_0}\right)+\sum_{j \in \cJ(B)}m_{j}\ln\left(\frac{Le}{m_j}\right)}{n}\right)^{\frac{\alpha+1}{\alpha+2}}
\ee

In addition,
$$
\sup_{B \in \cM(d_0)} \E\|V^{\frac{1}{2}}(\widehat{B}_{sgS}-B)\|^2_F \leq \frac{C_1}{\nu_{gS}(d_0)} \frac{d_0 \ln\left(\frac{de}{d_0}\right)+\sum_{j \in \cJ(B)}m_{j}\ln\left(\frac{Le}{m_j}\right)}{n}
$$
and
$$
\sup_{B \in \cM(d_0,{\bf m})}  \E||\widehat{B}_{sgS}-B||_{sgS} \leq
\frac{C_2}{\nu_{gS}(d_0)}~ \frac{d_0 \ln\left(\frac{de}{d_0}\right)+\sum_{j \in \cJ(B)}m_{j}\ln\left(\frac{Le}{m_j}\right)}{n}
$$
\end{corollary} 
Corollary \ref{cor:complexity_bound} shows that with a proper choice of tuning parameters, the bounds for misclassification excess risk for the global row-wise sparsity (\ref{eq:upper_noise}) are improved under the stronger double row-wise sparsity assumption.
The multinomial logistic sparse group Slope classifier $\widehat{\eta}_{sgS}$ is adaptive to $d_0, {\bf m}$ and $\alpha$.

Similar to global sparsity, there is a phase transition between small and large number of classes.
The numerator in the upper bounds  contains two terms. The first term 
$d_0\ln(de/d_0)$ corresponds again to the error of selecting a subset of $d_0$ nonzero rows out of $d$, while the second term 
$\sum_{j \in \cJ(B)}m_{j}\ln(Le/m_j)$
appears due to simultaneous estimation of $d_0$ $m_j$-sparse vectors from $\mathbb{R}^L$. 
Since $d_0\ln (Le) \leq \sum_{j \in \cJ(B)}m_{j}\ln(Le/m_j) \leq d_0 L$, the first term is always dominating for small number
of classes with $L \leq \ln(de/d_0)$, while the second term is the main one for large number of classes with $L \geq d/d_0$.

It also follows from Theorem \ref{th:sgs_low_noise_classifier} that, similar to the group Lasso, the multinomial logistic sparse group Lasso 
classifier with constant $\lambda=c_1\sqrt{\frac{\ln d}{n}}$ and $\kappa=c_2 \sqrt{\frac{\ln L}{n}}$ in (\ref{eq:sgS}) is sub-optimal up to the differences in the log-terms:
$$
\sup_{\eta^* \in \cC_L(d_0, {\bf m})}\cE(\widehat{\eta}_{sgL},\eta^*) 
\le  \left(\frac{C}{\nu_{gS}(d_0)}~\frac{d_0 \ln d +\ln L \cdot \sum_{j \in \cJ(B)}m_{j}}{n} \right)^{\frac{\alpha+1}{\alpha+2}}.
$$

Note that unlike global row-wise sparsity, 
interpretation of local (and, therefore, the double) row-wise sparsity assumption depends on the chosen constraint on $B$. Thus, a non-zero row of $B$ may be sparse (in terms of $l_0$-norm) under the symmetric constraint $\sum_{l=1}^L \bbeta_l={\bf 0}$ but not necessarily sparse under another possible constraint, e.g., $\bbeta_L={\bf 0}$ and vice versa. 

\ignore{
One can consider an unconstrained multinomial logistic sparse group Slope estimation in (\ref{eq:sgS}). Its solution evidently exists and is unique due to convexity.
For a given matrix $B$ in (\ref{eq:multinom}) consider a class of all matrices obtained by shifting each column $\bbeta_l$ of $B$ by
all possible vectors ${\bf c} \in \mathbb{R}^d$. As we have mentioned, all these matrices are equivalent in terms of the likelihood.
The unconstrained multinomial logistic sparse group Slope corresponds then to a general assumption of the existence of a double sparse matrix $B^*$ (not necessarily the one with the symmetric constraint)
within this equivalence class. Such a  relaxed assumption allows more flexibility but at a price of loosing possible interpretability of $B^*$. The following theorem shows that the upper bounds for misclassification excess risk remain of the same order: 
\begin{theorem} \label{th:sgs_low_noise_classifier_unconstrained}
\end{theorem} 
}

\ignore{
Moreover, it can be proved that the entry $c_j^*$ of the corresponding shift
vector ${\bf c}^*$ will lie between the mean and the median of the $j$-th row of $B$ (see Theorem 1 of Friedman et al., 2010).
}

\ignore{
The resulting (unique) solution
$\hat{B}^*_{sgS}$ 
}

\subsection{Low-rank sparsity} \label{subsec:rank}
So far we considered various types of row-wise sparsity of the regression coefficients matrix $B$. A more general approach is to assume the existence of some underlying hidden low-dimensional structure, where there is a smaller number of latent variables that define the outcome classes. The row-wise sparsity is a particular case of such a general case. The natural measure of such type of sparsity (sometimes called also {\em spectral} sparsity) is $rank(B)$.

Direct penalization of $rank(B)$ implies a non-convex optimization since $rank(B)$ is not a convex function although \citet{she2013reduced} proposed a computationally fast procedure for its solution for GLM. To convexify rank penalization note that $rank(B)=||\bgamma||_0$, where $\gamma_1,\ldots,\gamma_{\min(L-1,d)}$ are the singular values of $B$. Similar to
Lasso, we replace
$||\bgamma||_0$ by $||\bgamma||_1$ aka a nuclear norm $||B||_*$ or, Schatten $S_1$-norm. Nuclear penalties have been intensively studied in statistical and machine learning for multivariate regression and matrix completion \citep[e.g.,][]{bach2008consistency,candes2010matrix,bunea2011optimal,koltchinskii2011nuclearnorm,alquier2019estimation}. \citet{powers2018nuclear} considered nuclear penalization in multinomial logistic classification.
They developed numerical algorithms for its
solution but did not investigate theoretical properties of the resulting classifier.

We start from establishing a minimax lower bound for misclassification excess risk over a set of $L$-class linear classifiers with law rank coefficients matrices.
Let $\cM^*(r_0)=\{B \in \cM: rank(B) \leq r_0\}$ and
$\cC^*_L(r_0)=\{\eta(\bx)=\arg \max_{1 \leq l \leq L} \bbeta^T_l\bx: B \in \cM^*(r_0)\}$. 

\begin{theorem} \label{th:nuclear_lower}
Consider an agnostic multinomial regression model (\ref{eq:model})-(\ref{eq:multinom}) with $rank(B) \leq r_0$,
where $1 \leq r_0 \leq \min(L-1,d)$ and $r_0(L+d) \leq n$. Then,
\be \label{eq:nuclear_lower}
\inf_{\widetilde \eta}\sup_{\eta^* \in \cC^*_L(r_0),~\mathbb{P}_X}
\cE(\tilde{\eta}, \eta^{*}) \ge C \sqrt{\frac{r_0 ((L-1)+d)}{n}}
\ee
for some $C>0$.
\end{theorem}
The proof is given in the Appendix \ref{sec:appendix_lower}.

We now show that estimating $B$ by penalized maximum likelihood estimator with a nuclear penalty of the form $\lambda ||B||_*$ with a properly chosen tuning parameter $\lambda$ leads to a  linear classifier that achieves the lower bound (\ref{eq:nuclear_lower}) up to a multiplicative term depending on the marginal distribution $\mathbb{P}_X$ of $\bX$.
 
Define
\begin{equation} \label{eq:nuclearpen}
\widehat{B}_{nu}=\argmin_{\widetilde{B}} \left\{ \frac{1}{n}\sum_{i=1}^{n} \left(\ln\left(\sum_{l=1}^{L} \exp(\widetilde{\beta}_{l}^{\top} \bx_{i})\right) - \bbeta_{y_{i}}^{\top} \bx_i \right)+ \lambda ||B||_* \right\},
\end{equation}
with $\lambda>0$, and the corresponding
classifier $\widehat{\eta}_{nu}(\bx)=\argmax_{1 \leq l \leq L}\widehat{\bbeta}_{nu,l}^T \bx$. Similar to group Slope and group Lasso classifiers from Section \ref{subsec:global}, there is no need to impose an additional symmetric constraint $\widetilde{B} \in \cM$ in (\ref{eq:nuclearpen}) since centering rows to zero means can only decrease the nuclear norm of a matrix \citep{powers2018nuclear}.

Let $\tau_1(V) \geq \ldots \geq  \tau_d(V)$ be the ordered eigenvalues of the second moment matrix $V=\E_X(\bX \bX^T)$.
\ignore{
In the case of nuclear norm classifier, we simply define \(\nu = \tau_{d}(\bX) = \sigma_{min} (\E_{X} [\bX \bX^{\top}])\), and require the following assumption,
}

\begin{subassumption} \label{as:Dnu}
Assume that 
        \(\tau_{d}(V) > 0 \).
\end{subassumption}
\ignore{
Note that 
\[
\frac{\|V^{\frac{1}{2}} B \|^2_{F}}{\| B \|^2_{F}} \ge \tau_{d}(\bX).
\]
Thus, Assumption \ref{as:Dnu} with \(\nu = \tau_{d}(\bX)\) is in fact of the same form as \ref{as:Dgs} and \ref{as:Dsgs}.
}

\begin{theorem} \label{th:nuclear_upper}
Consider a multinomial regression model (\ref{eq:model})-(\ref{eq:multinom})
and the nuclear penalized classifier $\widehat{\eta}_{nu}(\bx)$ with
\begin{equation} \label{eq:lambdanu}
\lambda=C \sqrt{\tau_1(V)} ~\frac{\sqrt{L-1} + \sqrt{d}}{\sqrt{n}},
\end{equation}
where $C>0$ is specified in the proof.

Then, under Assumption \ref{as:Dnu}
\begin{equation} \label{eq:nuclear_upper1}
\sup_{\eta^* \in \cC^*_L(r_0)}\cE(\widehat{\eta}_{nu},\eta^*) 
\le \sqrt{C~ \frac{\tau_{1}(V)}{\tau_{d}(V)}~ \frac{r_0 ((L-1)+d)}{n}}~.
\end{equation}

Furthermore, under the additional low-noise Assumption \ref{as:B}, 
\begin{equation} \label{eq:nuclear_upper2}
\sup_{\eta^* \in \cC^*_L(r_0)}\cE(\widehat{\eta}_{nu},\eta^*) 
\le  \left(C~\frac{\tau_{1}(V)}{\tau_{d}(V)}~   \frac{r_0 ((L-1)+d)}{n}\right)^{\frac{\alpha+1}{\alpha+2}}~.
\end{equation}

In addition,
$$
\sup_{B \in \cM(d_0)} \E\|V^{\frac{1}{2}}(\widehat{B}_{nu}-B)\|^2_F \leq C_1~\frac{\tau_{1}(V)}{\tau_{d}(V)}~   \frac{r_0 ((L-1)+d)}{n}
$$
and
$$
\sup_{B \in \cM^*(r_0)} \E||\widehat{B}_{nu}-B||_* \leq 
C_2~\frac{\tau_{1}(V)}{\tau_{d}(V)}~   \frac{r_0 ((L-1)+d)}{n}
$$
\end{theorem}
The proof is given in the Appendix \ref{sec:appendix_upper}.

Similar upper bounds for the misclassification excess risk with the extra $\ln^{3/2}(n^{3/2} L)$-term can be derived from Corollary 10 of \citet{lei2019datadependent} using (\ref{eq:KL}) from Appendix \ref{sec:appendix_upper}.

\ignore{
\noindent
\newline
{\bf Remark} The upper bounds (\ref{eq:nuclear_upper1}) and  (\ref{eq:nuclear_upper2}) hold for a very general marginal distribution $\mathbb{P}_X$ of $\bX$ and can be improved  under additional assumptions on $\mathbb{P}_X$. Thus, for a sub-Gaussian $\bX$,
applying the results of \citet[Section 5.4.1]{vershynin2012introduction} in the proof it can be shown
that with $\lambda=C\sqrt{\tau_{1}(\bX)}~\frac{\sqrt{L-1} + \sqrt{d}}{\sqrt{n}}$
$$
\sup_{\eta^* \in \cC^*_L(r_0)}\cE(\widehat{\eta}_{nu},\eta^*) 
\le \left(C~ \frac{\tau_{1}(\bX)}{\tau_{d}(\bX)} \frac{r_0 \left((L-1)+d\right)}{n}\right)^{\frac{\alpha+1}{\alpha+2}}.
$$
}

Summarizing, up to a multiplicative constant depending on the eigenvalues of the second moment matrix of $\bX$, $\widehat{\eta}_{nu}(\bx)$ attains the minimax misclassification excess risk and is adaptive to the unknown low-rank sparsity of the regression coefficients matrix. 

\ignore{
The requirement $r_0((L-1)+d) \leq n$ necessary for consistent classification implies, in particular, that $d < n$ (although the number of linearly independent entries $d \times (L-1)$ of $B$ may be still larger than $n$ for sufficiently low $r_0$). For $d >n$ one should assume some extra assumptions on the sparsity of $B$ (e.g., row-wise sparsity) and consider combinations of row-wise and low-rank sparsities. The details are beyond the scope of the paper.
}

\section{Example} \label{sec:examples}
To illustrate the performance of the derived sparse multinomial logistic regression classifiers we applied them to the data set {\em Cancer sites} considered in  \citet{vincent2014sparse}.
It consists of bead-based expression data for $n=162$ microRNAs with $d=372$ features from $L=18$ classes of normal and cancer issue samples. The number of samples in each class ranges from 5 to 26. \citet{vincent2014sparse} used sparse group Lasso classifier for this data.

\ignore{
The {\em Amazon reviews} data set includes $d=10000$ various text features
extracted from $n=1500$ customer reviews on $L=50$ authors.
}

We compared the performance of sparse group Slope with $\lambda_j$'s and $\kappa_\ell$'s
of the form given in (\ref{eq:sgslope_parameters}), sparse group Lasso \citep[replicating][]{vincent2014sparse}, random forest and the well-known gradient boosting trees XGBoost classifiers on the above data set, where we developed the proximal gradient algorithm for solving sparse group Slope in (\ref{eq:sgS}) -- see Appendix \ref{sec:sgS_algorithm}.

To remove various technical variations, following Vincent and Hansen (2014), the data was first normalized by centering and scaling the rows of the design matrix, and then standardized by centering and scaling the columns.
We split the data into training (75\%) and test (25\%) sets.
The tuning parameters of all classification procedures were chosen by 10-fold cross-validation on the training set,  and  the misclassification errors of the resulting classifiers were measured on the test set. We repeated the process 10 times, randomly partitioning the data into train and test sets.

\ignore{
We developed the proximal gradient algorithm for solving (\ref{eq:sgS}) -- see Appendix \ref{sec:sgS_algorithm}.
The constant tuning parameters $c_1, c_2$ in (\ref{eq:sgslope_parameters}) for sparse group Slope were chosen by cross-validation. 

The code and the data are available at \url{https://github.com/tmrlvi/sparse.group.slope}.
}

Table \ref{tab:results} presents the average (over 10 random splits) misclassification errors for the test sets, the numbers of selected features (non-zero rows of the regression coefficients matrix $B$) and the overall numbers of non-zero coefficients in $B$. It shows that  both sparse multinomial logistic regression classifiers outperform their nonparametric counterparts for this data. Sparse group Slope yielded smaller misclassification errors than sparse group Lasso and, in addition, resulted in much sparser models. 

\begin{center}
\begin{table} 
\begin{tabular}{ |l||c|c|c|  }
\hline
 Classifier & Average misclass. error & \# features & \# non-zero coefficients\\
 \hline
 sparse group Slope  & 0.159 (0.019) & 60-67 & 186-271 \\
 sparse group Lasso & 0.165  (0.018) & 51-79  &  382-592 \\
 random forest & 0.209 (0.009) & - &   -   \\
 XGBoost & 0.250 (0.026) & - &   -   \\
 \hline
\end{tabular}
\caption{Average misclassification errors with their standard errors (in brackets) and feature selection for various classifiers.}
\label{tab:results}
\end{table}
\end{center}

\ignore{
\begin{center}
\begin{tabular}{ |l||c|c|c||c|c|  }
\hline
 Dataset & \(n\) & \(d\) & \(L\) & Sparse Group SLOPE & \vtop{\hbox{\strut Sparse Group LASSO}\hbox{\strut \smaller{\citep{vincent2014sparse}}}} \\
 \hline
 Cancer Sites   & \( 162 \) & \(217\) & \( 18 \) &  \textbf{14.12}  & 19.75 \\
 Amazon Reviews & \( 1500 \) & \(10,001\) & \( 50 \) &   ?  & ? \\
 \hline
\end{tabular}
\end{center}
}

\section{Concluding remarks} \label{sec:conclusion}
In this paper we discussed high-dimensional multiclass classification by sparse multinomial logistic regression.  Multiclass setup allows one to consider various types of sparsity associated with different assumptions on a matrix of regression coefficients. 
We proposed penalized MLE feature selection procedures with 
convex penalties capturing a specific type of sparsity at hand and showed
that the resulting classifiers are optimal in the minimax sense.
We presented the results for global row-wise,  double row-wise and low-rank sparsity scenarios but 
one can consider also other related types of sparsity, e.g., group-sparsity, when features may have a group structure, or class-dependent sparsity, where each class has its own sparse
subset of predictive features that implies column-wise sparsity, combinations of row-wise and low-rank sparsities, etc. The developed approach is general (see Appendix \ref{subsec:appendix_uppper_general} and Theorem \ref{th:generalization} there) although a specific type of a penalty should be properly chosen w.r.t. a particular type of sparsity at hand.

In this paper we assume that $\mathbb{P}_X$ has a bounded support. Using a slightly different techniques, the main results remain valid also for Gaussian design \citep[see][for binary classification]{bellec2018slope,alquier2019estimation}.

Even when the considered multinomial logistic regression model is misspecified and the Bayes
classifier $\eta^*$ is not linear, the misclassification excess risk can still be decomposed as
\begin{equation} \label{eq:decomp}
R(\widehat{\eta}_{\widehat B})-R(\eta^*)=\left(R(\widehat{\eta}_{\widehat B})-R(\eta^*_L)\right) +
\left(R(\eta^*_L)-R(\eta^*)\right),
\end{equation}
where $\eta^*_L=\arg \min_{\eta \in \cC_L} R(\eta)$ is the best possible (oracle) linear classifier.  
The results of the paper can then be applied to the first term in the RHS of (\ref{eq:decomp})
representing the estimation error, whereas the approximation error in the second term measures the ability of linear classifiers to perform as good as $\eta^*$.
Enriching the class of linear classifiers may improve
the approximation error but will increase the resulting estimation error in (\ref{eq:decomp}). 
In a way, it is similar to the variance/bias tradeoff in regression. 

\section*{Acknowledgments}
The work was supported by the Israel Science Foundation (ISF), Grants ISF-589/18 and ISF-1095/22.
The authors would like to thank Amir Beck, Guillaume Lecu\'e and the anonymous referees for helpful comments.

\appendix

\section{Proofs of the upper bounds (Theorems \ref{th:gs_low_noise_classifier}, \ref{th:sgs_low_noise_classifier} and  \ref{th:nuclear_upper})} \label{sec:appendix_upper}

Throughout the proofs we use various generic positive constants, not necessarily the same each
time they are used even within a single equation.

Throughout the proofs let $|{\bf a}|_2$ be the Euclidean norm of a vector ${\bf a}$, $||A||_2$ and $||A||_F$ respectively the operator/spectral and Frobenius norms of a matrix $A$. 
The Frobenius inner product of two matrices $A_1$ and $A_2$ is $\langle A_1, A_2\rangle =tr(A_1^\top A_2)$.
Denote $||g(\bx)||_{L_2}$ for the $L_2$-norm of a function $g$ and $||g(\bx)||_{L_2(\mathbb{P}_X)}
=(\int_{\cX} g(\bx)^2 d\mathbb{P}_X(\bx))^{1/2}$ for the $L_2$-norm of $g$  w.r.t. the measure $\mathbb{P}_X$. Recall that \(V= \E[\bX\bX^{\top}]\).
\ignore{
In addition, let $||A||_{S_p},\; p \geq 1$ be a Schatten $p$-norm of a matrix $A$. In particular, $||A||_{S_1}=||A||_*$ (nuclear norm), $||A||_{S_2}=||A||_F$ and $||A||_{S_\infty}=||A||_2$.
}
\subsection{Upper bounds for misclassification excess risk for a general penalized MLE plug-in linear classifier} \label{subsec:appendix_uppper_general}
  Consider first a generic setup. 
  Let \(\cM = \{ B \in \R^{d \times L} : B{\bf 1}={\bf 0} \}\) be the set of regression matrices satisfying the symmetric constraint and
  \(\cM_0 \subseteq \cM\) be its subset of sparse matrices, where the notion of sparsity depends on the particular problem at hand. Let $B \in \cM_0$ and consider a penalized MLE estimator $\widehat{B}$ of the form
\begin{equation} \label{eq:generalPMLE}
\widehat{B}=\arg \min_{\widetilde{B} \in \cM} \left\{
-l(\widetilde{B})+||\widetilde{B}|| \right\},
\end{equation}
where the regularized matrix norm $||\cdot||$ induces the given type of sparsity, and the corresponding plug-in linear classifier 
$$
\widehat{\eta}_{\widehat B}(\bx)=\arg \max_{1 \leq l \leq L} \widehat{\bbeta}_l^\top \bx.
$$

The Kullback-Leibler divergence between two multinomial distributions with probabilities vectors ${\bf p}_1$ and ${\bf p}_2$ is $KL({\bf p}_1,{\bf p}_2)=\sum_{l=1}^L p_{1l}\ln \left(\frac{p_{1l}}{p_{2l}}\right)$.
Let $f_B(\bx,y)$ be the joint distribution of $(\bX,Y)$, i.e., $df_B(\bx,y)=\prod_{l=1}^L p_l(\bx)^{\xi_l} d\mathbb{P}_X(\bx)$, where $p_l(\bx)$ are given in (\ref{eq:multinom}).
For two given regression coefficients matrices $B_1$ and $B_2$ the Kullback-Leibler divergence between the distributions $f_{B_1}$ and
$f_{B_2}$ is then $d_{KL}(f_{B_1},f_{B_2})=\int KL({\bf p}_1(\bx),{\bf p}_2(\bx)) d\mathbb{P}_X(\bx)$.
We exploit the well-known result \citep[e.g., ][]{pires2016multiclass,abramovich2021multiclass} that relates the misclassification excess risk $\cE(\hat{\eta}, \eta^*)$ and the Kullback-Leibler risk $\E d_{KL}(f_B,f_{\widehat B})$ under the low-noise Assumption \ref{as:B}:
\be \label{eq:KL}
\cE(\widehat{\eta}_{\widehat M},\eta^*) \leq C 
\left(\E d_{KL}(f_B,f_{\widehat B}\right)^{\frac{\alpha+1}{\alpha+2}}.
\ee

We  now extend the results of \cite{alquier2019estimation} for univariate response to multivariate (multinomial) $\bY$ to bound the Kullback-Leibler risk $\E d^2_{KL}(f_B,f_{\widehat B})$. 
Define $\theta_l(\bx)=\bbeta_l^\top \bx,\;l=1,\ldots,L$, where due to the symmetric constraint, $\sum_{l=1}^L \theta_l(\bx)=0$. It is easy to verify that in terms of $\theta_l$'s, the multinomial log-likelihood is Lipschitz w.r.t. the $l_2$-norm. Furthermore, for $\mathbb{P}_X$ with a bounded support, $|\bbeta_l^\top \bx| \leq C$ and  $d_{KL}(\cdot,\cdot)$ is strongly convex \citep{abramovich2021multiclass}: for any two matrices $B_1$ and $B_2$ satisfying the symmetric constraint, $d_{KL}(f_{B_1},f_{B_2}) \geq C \sum_{l=1}^L  ||\theta_{1l}(\bx)-\theta_{2l}(\bx)||^2_{L_2(\mathbb{P}_X)}$  (the multivariate analogue of Bernstein condition in terminology of \citep{alquier2019estimation}.
These two conditions allow us to adopt the general approach of \citet{alquier2019estimation} to bound $\E d^2_{KL}(f_B,f_{\widehat B})$.

Define the following quantities along the lines of \citet{alquier2019estimation}. 
Let $\cb_{||\cdot||}=\{B \in \cM: ||B|| \leq 1\}$ be the unit ball of matrices satisfying the symmetric constraint w.r.t. $||\cdot||$-norm in (\ref{eq:generalPMLE}).
Let $\widehat{Rad}(\cb_{||\cdot||})$ be the empirical (multivariate) Rademacher complexity of $\cb_{||\cdot||}$, namely,
\begin{equation} \nonumber
\begin{split}
\widehat{Rad}(\cb_{||\cdot||})&=\E_\Sigma\left\{\frac{1}{\sqrt n}\sup_{B 
\in \cb} \sum_{i=1}^n \sum_{l=1}^{L} \sigma_{il} \bbeta_l^T \bX_i \Big|\bX_1=\bx_1,\ldots,\bX_n=\bx_n \right\}\\
&=\E_\Sigma \left\{\frac{1}{\sqrt n}\sup_{B 
\in \cb} tr(\Sigma B^T X^T)\right\},
\end{split}
\end{equation}
where the elements $\sigma_{il}$'s of $\Sigma \in \mathbb{R}^{n \times L}$ are i.i.d.
Rademacher random variables with $P(\sigma_{il}=1)=P(\sigma_{il}=-1)=1/2$, and 
$$
Rad(\cb_{||\cdot||})=\E_X\left\{\widehat{Rad}(\cb_{||\cdot||})\right\}
$$
be the Rademacher complexity of $\cb$.

Define a {\em complexity function}
$$ 
\complexity(\rho) = \sqrt{\frac{C_0 Rad(\cb_{||\cdot||}) \rho}{2 R^{2} \sqrt{n}}},\;\;\;\rho>0,
$$
where the exact value of $C_0>0$ is specified in \cite{alquier2019estimation}.

Let $\mathcal{T}(\rho) = \{B^\prime \in \cM: ||B^\prime||= \rho,~ ||V^{\frac{1}{2}}B^\prime||^2_F \le r^2(2\rho)\}$.  For a given matrix $B \in \cM_0$ define $\Gamma_{B}(\rho) = \bigcup_{B' : || B' - B|| < \frac{\rho}{20}} \partial ||\cdot ||(B')$, where the subdifferential $\partial ||\cdot ||(B')=\{G \in \cM: ||B^\prime+B^{\prime \prime}||- ||B^\prime|| \geq \langle B^{\prime\prime}, G \rangle,\; \forall B^{\prime\prime} \in \cM\}$. The {\em sparsity parameter} is 
   $$
    \sparsity(\rho) =\inf_{B' \in \mathcal{T}(\rho)} \sup_{G \in \Gamma_{B}(\rho)} \langle B^\prime, G \rangle.
    $$
Finally, let $\rho^*$ be any solution of the sparsity inequality
\begin{align}\label{eq:sparsity_equtation}
\sparsity(\rho^*) \geq \frac{4}{5}\rho^*
\end{align}

The quantity $\rho^*$ depends on a particular norm in (\ref{eq:generalPMLE}) and the second moment matrix $V$, and plays a key role in establishing the upper bound for  $Ed_{KL}(f_{B_1},f_{B_2})$.  

We have the following generic theorem:
\begin{theorem} \label{th:generalization}
	Let $\widehat{B}$ be the solution of (\ref{eq:generalPMLE}). Assume that there exists $\rho^*$
	such that $\sparsity(\rho^*) \geq \frac{4}{5}\rho^*$ and $Rad(\cB_{||\cdot||}) \le \frac{7}{720} \sqrt{n}$. 
	Then, 
	\begin{equation} \label{eq:general_upper}
    \E d_{KL}(f_B, f_{\widehat B}) \le C \rho^*,
	\end{equation}
	for some $C>0$. 
	
	In addition,
	$$
    \E\| V^{\frac{1}{2}}(\widehat{B}-B)\|^2_F
	\le C_1 \rho^*
	$$
	and
	$$
    \E\|\widehat{B} - B\| \le C_2 \rho^*
	$$
	for some $C_1, C_2>0$.
\end{theorem}

\ignore{
\begin{assumption} \label{as:D} Consider a subset of sparse matrices $\cM_0 \subseteq \cM$, where the notion of sparsity depends on a particular problem at hand, and let $\nu(\cM_0)$ be the corresponding minimal sparse eigenvalue of $V$:
$$
\nu(\cM_0) =\inf_{A \in \cM_0: A \neq 0_{d \times L}} \frac{||V^{1/2}A||^2_F}{||A||^2_F}
$$
Assume that $\nu(\cM_0) >0$.
\end{assumption}
As we have mentioned, such type of assumptions are always required for convex relaxation methods. 
}
 
Theorem \ref{th:generalization} is an extension of Theorem 2.2 (or more general Theorem 9.2)
of \citet{alquier2019estimation} for multivariate response and anisotropic design. 
Its proof repeats the proof of Lemma 1 in \citet{abramovich2021multiclass} with the particular group Slope norm considered there replaced by a general norm $||\cdot||$.

\begin{remark} \label{rem:general}  
{\rm
 In fact,  from the definition of the sparsity parameter $\Delta(\rho)$ and $\rho^*$ it follows that Theorem \ref{th:generalization} holds even if the true regression matrix $B$ is only ``approximately sparse'' in the sense that there exists a sparse matrix $B' \in \cM_0$ such that $||B-B'|| \leq \rho^*/20$ \citep[see also][]{alquier2019estimation}.}
\end{remark}

We will now apply the general upper bound (\ref{eq:general_upper}) for the group Slope, sparse group Slope and nuclear norms to complete the proofs of Theorems \ref{th:gs_low_noise_classifier}, \ref{th:sgs_low_noise_classifier} and  \ref{th:nuclear_upper} by finding the corresponding  $Rad(\cB_{||\cdot||})$ and $\rho^*$.

\subsection{Proof of Theorem \ref{th:gs_low_noise_classifier}}
 The proof of Theorem \ref{th:gs_low_noise_classifier} is somewhat different from that of Theorem 4 of \citet{abramovich2021multiclass} for isotropic $\bX$.
\ignore{
The proof of Theorem \ref{th:gs_low_noise_classifier} differents from the proof by \citet{abramovich2021multiclass} in the calculation of the sparsity parameter. There, instead of extending Lemma 4.3 of \cite{lecue2018regularization}, we extend Lemma 0.2 of their supplementary material, while noting that assumption \ref{as:Dgs} is the matrix equivalent to their second part of Assumption 0.1. We lay out the argument here for completeness.
}

For given $\lambda_1 \geq \ldots \ge \lambda_d$ consider the group Slope norm
$||B||_\lambda=\sum_{j=1}^d \lambda_j |B|_{(j)}$. Let \(B \in \cM(d_0)\) with a set of zero rows \(\cJ(B)\)  and \(B^{\prime} \in \R^{d\times L}\) such that \(\|B^{\prime} - B\|_{\lambda} = \rho^{*}\) and 
\(\| V^{\frac{1}{2}}(B'-B)\|^2_F
\le \frac{C_0 Rad(\cB_{\lambda}) \rho^{*}}{ R^{2} \sqrt{n}}\), where $\rho^*$ will be defined later .

Let $\mathcal{G}$ be a set all of matrices of the form  $\sum_{j \in \cJ(B)} \lambda_{\pi(j)} {\bf e}_{j} \frac{B_{j\cdot}}{|B_{j\cdot}|_2} + \sum_{j \in \cJ^c(B)} \lambda_{\pi(j)} {\bf e}_{j} {\bf v}_j^{\top}$, where  $\pi=(\pi(1),\ldots,\pi(d))$ is a permutation of $\{1,\ldots,d\}$ and ${\bf v}_j$'s are unit vectors in \(\R^{L}\), and note that
$
\| B \|_{\lambda} = \max_{G \in \mathcal{G}} \langle B, G\rangle.
$ 
\ignore{
where 
\[
\mathcal{G} = \left\{ 
\begin{array}{c}
\sum_{j \in \cJ(B)} \lambda_{\pi(j)} \frac{B_{j\cdot}}{\|B_{j\cdot}\|} e_{j}^{\top} + \sum_{j \in \cJ^c(B)} \lambda_{\pi(j)} v_j e_{j}^{\top} \\
\;\text{for all permutations}\; \pi=(\pi(1),\ldots,\pi(d)) \;\text{of}\; (1,\ldots,d) \\ \text{and unit vectors}\; v_{j} 
\end{array}
\right\}.
\]
}

In particular,
$
\partial \| \cdot \|_{\lambda}(B) \supseteq \argmax_{G \in \mathcal{G}} \langle B, G \rangle.
$
Hence, we can find a permutation of \(\{\lambda_{j}\}_{j=d_0+1}^{d}\) such that the corresponding \(G \in \argmax_{G \in \mathcal{G}} \langle B, G \rangle \subseteq \partial \| \cdot \|_{\lambda}(B)\) and  \(\sum_{j \in \cJ^{c}(B)} G_{\cdot j}^{\top} (B^{\prime} - B)_{j\cdot } \ge \sum_{j=d_0 + 1}^{d} \lambda_j |B^{\prime} - B|_{(j)} \). Then,
\begin{equation} \label{eq:rho*1}
\begin{split}
\left\langle G,B^{\prime} - B\right\rangle &= \sum_{j \in\cJ(B) } G_{\cdot j}^{\top} (B^{\prime} - B)_{j\cdot } + \sum_{j \in \cJ^{c}(B)} G_{\cdot j}^{\top} (B^{\prime} - B)_{j\cdot } \\
&\ge \sum_{j \in \cJ^{c}(B)} G_{\cdot j}^{\top} (B^{\prime} - B)_{j\cdot } - \sum_{j=1}^{d_0} \lambda_{j} | B^{\prime} - B |_{(j)} \ge \\
& \ge \sum_{j=1}^{d} \lambda_{j} | B^{\prime} - B |_{(j)} - 2\sum_{j = 1}^{d_0} \lambda_{j} | B^{\prime} - B |_{(j)} \\
&= \rho^{*} - 2 \sum_{j=1}^{d_0} \lambda_{j} |B^{\prime} - B |_{(j)}.
\end{split}
\end{equation}
By Assumption \ref{as:Dgs},
\begin{equation} \label{eq:rho*2}
\frac{1}{\nu_{gS}(d_0)}\| V^{\frac{1}{2}} (B^{\prime} - B)\|^2_{F} \ge \|\Pi_{d_0}(B^{\prime}-B)\|^2_F = \sum_{j=1}^{d_0} | B^{\prime} - B |_{(j)}^{2}~.
\end{equation}
For any $1 \leq j \leq d_0$ we also have
$$
\sum_{j^{\prime}=1}^{d_0}|B^{\prime}-B|_{(j^{\prime})}^2 \geq \sum_{j^{\prime}=1}^j|B^{\prime}-B|_{(j^{\prime})}^2 \geq j |B^{\prime}-B|_{(j)}^2
$$
and, therefore,
\(|B^{\prime} - B |_{(j)} \le \sqrt{\sum_{j^{\prime}=1}^{d_0} | B^{\prime} - B |_{(j^{\prime})}^{2}}/ \sqrt{j}\). 

Taking
$$
\rho^{*} = \frac{100 C_0 C}{\nu_{gS}(d_0)} \frac{Rad(\cB_{\lambda})\left(\sum_{j=1}^{d_0} \lambda_j/\sqrt{j}\right)^2}{\sqrt{n}}
$$
(\ref{eq:rho*2}) implies 
\[
\begin{split}
\sum_{j=1}^{d_0} \lambda_{j} | B^{\prime} - B|_{(j)} &\le \frac{1}{\sqrt{\nu_{gS}(d_0)}}
\left(\sum_{j=1}^{d_0}\frac{\lambda_{j}}{\sqrt{j}}\right)
\|V^{\frac{1}{2}}(B'-B)\|_F
\\
& \le \frac{1}{\sqrt{\nu_{gS}(d_0)}}\left(\sum_{j=1}^{d_0}\frac{\lambda_{j}}{\sqrt{j}}\right) \sqrt{\frac{C_0 Rad(\cB_{\lambda} ) \rho^{*}}{ R^{2} \sqrt{n}}} \\
& \le \frac{1}{10}\rho^*.
\end{split}
\]
Thus, combining with (\ref{eq:rho*1}) 
\[
\left\langle G,B^{\prime} - B\right\rangle \ge \frac{4}{5} \rho^*
\]
for every \(B^{\prime}-B \in \mathcal{T}(\rho^*)\) and, therefore,
\[
\Delta(\rho^{*}) \ge \frac{4}{5} \rho^{*}.
\]

Furthermore, by Lemma 2 of \citet{abramovich2021multiclass},
$$
Rad(\cB_{\lambda}) \leq C \max_{1 \leq j \leq d}\frac{\sqrt{L+\ln(d/j)}}{\lambda_j}
$$
for some $C>0$. Hence, for $\lambda_j$ satisfying (\ref{eq:weightsgs}), $Rad(\cB_{\lambda}) \leq \frac{7}{720}\sqrt{n}$ and we can apply Theorem \ref{th:generalization} to complete the proof.

\subsection{Proof of Theorem \ref{th:sgs_low_noise_classifier}} \label{subsec:appendix_upper}
For given $\lambda_1 \geq \ldots \geq \lambda_d>0$ and $\kappa_1 \geq \ldots \geq \kappa_L>0$, consider the sparse group Slope norm $||B||_{\kappa,\lambda}=\sum_{j=1}^d \lambda_j |B|_{(j)}+
\sum_{j=1}^d \sum_{l=1}^L \kappa_l |B|_{j(l)}$, where $|B|_{(1)} \geq \ldots \geq |B|_{(d)}$ are the descendingly ordered $l_2$-norms of the rows of $B$ and $|B|_{j(1)} \geq \ldots \geq B_{j(L)}$ are descendingly ordered absolute values of entries of its rows. Let $\cB_{\kappa,\lambda}$ be the unit ball of matrices w.r.t. this norm. 

\begin{lemma} \label{lem:rho_star_sgs}
Let $B \in \cM(d_0,{\bf m})$. Under Assumption \ref{as:Dgs}, define
\begin{equation} \label{eq:rho_star_sgs}
\rho^*=\frac{100 C_0 C}{\nu_{gS}(d_0)}~ \frac{Rad(\cb_{\kappa,\lambda})\left(\sum_{j=1}^{d_0}\frac{\lambda_j}{\sqrt{j}}+\sqrt{\sum_{j=1}^{d_0} \left(\sum_{l=1}^{m_j}\frac{\kappa_l}{\sqrt{l}} \right)^2} \right)^2}{\sqrt{n}}.
\end{equation}
Then, $\rho^*$ satisfies the sparsity inequality (\ref{eq:sparsity_equtation}), i.e., $\sparsity(\rho^*) \geq \frac{4}{5}\rho^*$.
\end{lemma}

To apply Theorem \ref{th:generalization} to complete the proof, we need also to show that \(Rad(\cB_{\kappa,\lambda}) \le \frac{7}{720}\sqrt{n}\):
\begin{lemma} \label{lem:rademacher_sgs} Let \(\kappa_{L} \ge \sqrt{\frac{\pi}{2}}\frac{2880}{7} \frac{1}{\sqrt{n}}\).
Then, 
	$$
	Rad(\cB_{\kappa,\lambda}) \le \frac{7}{1440} \sqrt{n} + C_{0} \sqrt{\frac{\pi}{2}}\max_{1\le j \le d} \frac{\sqrt{2 \sum_{j=1}^{L}\frac{1}{l}\left(\frac{Le}{l}\right)^{l}e^{-C^2 n l \kappa_{l}^{2}}+2\log\left(\frac{de}{j}\right)}}{\lambda_{j}},
	$$
\end{lemma}
where $C = \sqrt{\frac{2}{\pi}}\frac{7}{2880}$ and $C_0>0$ is given in the proof.

In particular, for $\lambda_j$'s and $\kappa_l$'s satisfying (\ref{eq:weightssgs}), \(Rad(\cB_{\kappa, \lambda}) \le \frac{7}{720} \sqrt{n}\).

\ignore{
Clearly, $e^{-C^2 nl\kappa_{l}^{2}} = \left(\frac{Le}{l}\right)^{-2l}$.
Therefore, 
$$
\sum_{l=1}^{L}\frac{1}{l}\left(\frac{Le}{l}\right)^{l}e^{-C^2 nl\kappa_{l}^{2}}=\sum_{l=1}^{L}\frac{1}{l}\left(\frac{l}{Le}\right)^{l}\le\sum_{l=1}^{L}\frac{1}{l}\frac{l}{Le}=\frac{1}{e}
$$
Finally, 
\[
\max_{j}\frac{\sqrt{2\sum_{l=1}^{L}\frac{1}{l}\left(\frac{Le}{l}\right)^{l}e^{-2nl\kappa_{l}^{4}}+\log\left(\frac{de}{j}\right)}}{\lambda_{j}}\le\max_{j} \frac{7}{1440 \sqrt{2} C_0} \sqrt{\frac{2n}{e\log\left(de/j\right)}+n}\le \frac{7}{1440 C_0} \sqrt{n}.
\]
Thus, by Lemma \ref{lem:rademacher_sgs}, $Rad(\cB_{\kappa,\lambda}) \leq \frac{7}{720}\sqrt{n}$ and substituting $\rho^*$ from (\ref{eq:rho_star_sgs}) into (\ref{eq:general_upper}) completes the proof.
}

\subsection{Proof of Theorem \ref{th:nuclear_upper}}
Let $||B||_\lambda=\lambda ||B||_*$ and $\cB_\lambda$ the corresponding unit ball.
Define 
	\[ 
	\rho^* = 100 \lambda \frac{C_0 r_0 Rad(\cb_\lambda)}{2 R^{2} \tau_{d}(V) \sqrt{n}}.
	\]
	Extending Lemma 4.4 of \citet{lecue2018regularization} for the anisotropic case by using \(\| B\|_* < \frac{1}{\sqrt{\tau_{d}(V)}}\| V^{\frac{1}{2}} B\|_* \), we have $\sparsity(\rho^*) \ge \frac{4}{5} \rho^*$.
	
	To apply Theorem \ref{th:generalization} we need to show that for $\lambda$ in (\ref{eq:lambdanu}),
$Rad(\cB_\lambda) \le \frac{7}{720} \sqrt{n}$:
\begin{lemma} \label{lem:nu_rademacher} 
	$$
		Rad(\cB_\lambda) \le  C_0 \sqrt{\tau_{1}(V)}~ \frac{\sqrt{L-1} + \sqrt{d}}{\lambda}
	$$
	for some $C_0>0$.
\end{lemma}
Thus, taking \(C = \frac{720 C_0}{7}\), the choice of $\lambda = \sqrt{\tau_1(V)}~\frac{\left(\sqrt{L-1} + \sqrt{d}\right)}{\sqrt n}$ implies $Rad(\cB_{\lambda}) \le \frac{7}{720} \sqrt{n}$.

\section{Proofs of lemmas}

\subsection{Proof of Lemma \ref{lem:rho_star_sgs}}
We use the arguments similar to those in the proof of Theorem \ref{th:gs_low_noise_classifier}.

Let \(\cJ\) be the set of indices of non-zero rows of \(B\) and \(\cl_j\) be the set of indices
of non-zero entries of the \(j\)-th row for \(j \in \cJ\). Obviously, $|\cJ|=d_0$ and $|\cl_j|=m_j$. Consider a matrix \(B^{\prime}\) such that \(\| B^{\prime}-B \|_{\kappa, \lambda}=\rho^*\) and 
\(\|V^{\frac{1}{2}}(B^{\prime}-B)\|^2_F
\le r^2(2\rho^*) = \frac{C_0 Rad(\cB) C}{\sqrt{n}} \rho^* \).

We can decompose \(\| B\|_{\kappa, \lambda}\) into two additive components: \(\|B\|_{\kappa, 0}= \sum_{j=1}^{d}\sum_{l=1}^{L}\kappa_l|B|_{j(l)}\)
and \(\|B\|_{0, \lambda}=\sum_{j=1}^{d}\lambda_j |B|_{(j)}\). 

\ignore{Following the characterization by \citet{lecue2018regularization} of the extreme points of the Slope dual ball by  (p. 629), we can take \(G \in \cB_{\| \cdot \|^{*}_{\kappa,0}}, H \in \cB_{\| \cdot \|^{*}_{0,\lambda}}\) that satisfy
\[
\begin{aligned}
    B_{jl} \neq 0 \iff G_{jl} = \kappa_{\pi_{j}(l)} \sign(B_{jl}),
     \\ 
    | B_{j\cdot} |_2 \neq 0 \iff H_{jl} = \lambda_{\psi(j)} \frac{B_{jl}}{|B_{j \cdot} |_2},
\end{aligned}
\]
where \(\pi_{j}(l)\) denotes the order index of \(|B_{jl}|\) in \(\{|B_{jl^{\prime}}|\}_{l^{\prime}=1}^{L}\) when ordered, and \(\psi(j)\) denote the index of \(\|B_{j\cdot}\|\) in \(\{\|B_{j^{\prime}\cdot}\|\}_{j^{\prime}=1}^{d}\) when ordered. 
}

Define two matrices $G, H \in \mathbb{R}^{d \times L}$ as follows. For every $j \in {\cal J}$ let
$\pi_j(1),\ldots,\pi_j(m_j)$ be the indices of descendingly ordered nonzero entries $|B|_{j(l)}$'s and set
$G_{j\pi_j(l)}=\kappa_{\pi_{j}(l)} \sign(B_{j\pi_j(l)})$. Similarly, let 
$\tilde{\pi}(1),\ldots,\tilde{\pi}(d_0)$ be the indices of descendingly ordered 
Euclidean norms $|B|_{(j)}$ of $d_0$ nonzero rows of $B$ and set
$H_{jl} = \lambda_{\tilde{\pi}(j)} \frac{B_{\tilde{\pi}(j)l}}{|B_{\tilde{\pi}(j) \cdot} |_2}$. The entries of $G$ and $H$ corresponding to zero entries of $B$ will be defined later. 

By construction, \(tr(G^\top B)= \| B \|_{\kappa,0}\) and \(tr(H^\top B) = \| B \|_{0,\lambda}\), while for any \(B^{\prime}\), \(tr(G^\top B^{\prime}) \le \| B^{\prime} \|_{\kappa,0}\) and \(tr(H^\top B^{\prime}) \le \| B^{\prime} \|_{0,\lambda}\). Thus, \(G\) and \(H\) are in \(\partial \|\cdot\|_{\kappa,0} (B)\) and \(\partial \|\cdot\|_{0,\lambda}(B)\) respectively. 

We have 
\[\sum_{j = 1}^{d}\sum_{l\in \cl_{j}}G_{jl} \left| B_{jl}^{\prime}-B_{jl}\right| \le \sum_{j=1}^{d}\sum_{l=1}^{m_{j}}\kappa_{l}\left|B^{\prime}-B \right|_{j\left(l\right)}
\]
and
\[\sum_{j\in \cJ}\left| \sum_{l=1}^{L}H_{jl}\left(B_{jl}^{\prime}-B_{jl}\right) \right| \le \sum_{j=1}^{d_{0}}\lambda_{j}\left| B^{\prime}-B \right|_{(j)}.\]
Hence,
\begin{equation}
\begin{split} \label{eq:work1}    
tr\left(G^\top\left(B^{\prime}-B\right)\right) =& \sum_{j = 1}^{d}\sum_{l\in \cl_{j}}G_{jl}\left(B_{jl}^{\prime}-B_{jl}\right) + \sum_{j=1}^{d}\sum_{l\in \cl_{j}^{C}}G_{jl}\left(B_{jl}^{\prime}-B_{jl}\right) \\
\ge& \sum_{j=1}^{d}\sum_{l\in \cl_{j}^{C}}G_{jl}\left(B_{jl}^{\prime}-B_{jl} \right) - \sum_{j=1}^{d}\sum_{l=1}^{m_{j}}\kappa_{l}\left|B^{\prime}-B \right|_{j\left(l\right)},
\end{split}
\end{equation}
and
\begin{equation}
\begin{split} \label{eq:work2}    
tr\left(H^\top\left(B^{\prime}-B\right)\right) =& \sum_{j\in \cJ}\sum_{l=1}^{L}H_{jl}\left(B_{jl}^{\prime}-B_{jl}\right)+\sum_{j\in \cJ^{C}}\sum_{l=1}^{L}H_{jl}\left(B_{jl}^{\prime}-B_{jl}\right) \\
\ge& \sum_{j\in \cJ^{C}}\sum_{l=1}^{L}H_{jl}\left(B_{jl}^{\prime}-B_{jl}\right)-\sum_{j=1}^{d_{0}}\lambda_{j}\left| B^{\prime}-B \right|_{(j)}
\end{split}
\end{equation}
To bound the first terms of the RHSs in (\ref{eq:work1}) and (\ref{eq:work2})  from below for a given \(B^{\prime}\) complete the entries of \(G\) and \(H\) corresponding to zero entries of $B$ in such a way that
\[
\begin{aligned}
\sum_{j=1}^{d}\sum_{l\in \cl_{j}^{C}}G_{lj}\left(B_{lj}^{\prime}-B_{lj}\right) & \ge\sum_{j=1}^{d}\sum_{l=m_{j}+1}^{L}\kappa_{l}\left|B^{\prime}-B\right|_{j\left(l\right)},
\end{aligned}
\]
and
\[
\begin{aligned}
\sum_{j\in \cJ^{C}}\sum_{l=1}^{L}H_{lj}\left(B_{lj}^{\prime}-B_{lj}\right) & \ge\sum_{j=d_{0}+1}^{d}\lambda_{j}\left| B^{\prime}-B\right| _{(j)}.
\end{aligned}
\]

Thus, 
\[
 \begin{aligned}
tr\left(G^\top\left(B^{\prime}-B\right)\right) & \ge\sum_{j=1}^{d}\sum_{l=1}^{L}\kappa_{l}\left|B^{\prime}-B \right|_{j\left(l\right)}-2\sum_{j=1}^{d}\sum_{l=1}^{m_{j}}\kappa_{l}\left|B^{\prime}-B \right|_{j\left(l\right)} \\
 & =\| B^{\prime}-B \|_{\kappa,0} - 2\sum_{j=1}^{d}\sum_{l=1}^{m_{j}}\kappa_{l}\left|B^{\prime}-B \right|_{j\left(l\right)},
\end{aligned}
\]
and
\[
\begin{aligned}
tr\left(H^\top\left(B^{\prime}-B\right)\right) & \ge\sum_{j=1}^{d}\lambda_{j}\left| B^{\prime}-B\right|_{(j)}-2\sum_{j=1}^{d_{0}}\lambda_{j}\left| B^{\prime}-B\right|_{(j)}\\
 & = \| B^{\prime}-B \|_{0, \lambda} -2\sum_{j=1}^{d_{0}}\lambda_{j}\left| B^{\prime}-B\right|_{(j)}.
 \end{aligned}
 \]
Consider $Z=G+H$.  Evidently, \(Z \in \partial \|\cdot\|_{\kappa, \lambda}(B) \) and
\[
\begin{aligned}
tr \left(Z^\top\left(B^{\prime}-B\right)\right) & \ge \| B^{\prime}-B \|_{\kappa, \lambda} -2\sum_{j=1}^{d_{0}}\lambda_{j}\left| B^{\prime}-B\right| _{(j)}-2\sum_{j=1}^{d}\sum_{l=1}^{m_{j}}\kappa_{l}\left|B^{\prime}-B \right|_{j\left(l\right)}\\
 & = \rho^*-2\sum_{j=1}^{d_{0}}\lambda_{j}\left| B^{\prime}-B\right|_{(j)}-2\sum_{j=1}^{d}\sum_{l=1}^{m_{j}}\kappa_{l}\left|B^{\prime}-B \right|_{j\left(l\right)}.
\end{aligned}
\]

By Assumption \ref{as:Dgs},
\begin{equation*}
\frac{1}{\nu_{gS}(d_0)}\| V^{\frac{1}{2}} (B^{\prime} - B)\|^2_{F} \ge \|\Pi_{d_0}(B^{\prime}-B)\|^2_F = \sum_{j=1}^{d_0} | B^{\prime} - B |_{(j)}^{2} \ge \sum_{j \in \cJ} \sum_{l=1}^{m_j} | B^{\prime} - B |_{{j(l)}}^{2} ~.
\end{equation*}
\ignore{
Consider  a projection \(Proj_{\cM(d_0)}\) such that  \(\left| B^{\prime}-B\right|_{(j)} \le\frac{\left\| Proj_{\cM(d_0)}(B^{\prime}-B)\right\|_{F}}{\sqrt{j}}\) for all \(1\le j \le d_0\),
and  \(\left|B^{\prime}-B \right|_{j\left(l\right)}\le\frac{\left| Proj_{\cM(d_0)}(B^{\prime}-B) \right|_{2j}}{\sqrt{l}}\) for all  \(1 \le l \le m_j \).
}
Since \(|B^{\prime} - B |_{(j)} \le \sqrt{\sum_{j^{\prime}=1}^{d_0} | B^{\prime} - B |_{(j^{\prime})}^{2}}/ \sqrt{j}\),
\[
\sum_{j=1}^{d_{0}}\lambda_{j}\left| B^{\prime}-B\right| _{(j)} \le  \frac{1}{\nu_{gS}(d_0)}\sum_{j=1}^{d_{0}}\frac{\lambda_{j}}{\sqrt{j}} \| V^{1/2}(B^{\prime}-B) \|_F.
\]
On the other hand, by the Cauchy-Schwartz inequality,
\[ \begin{aligned}
\sum_{j=1}^{d}\sum_{l=1}^{m_{j}}\kappa_{l}\left|B^{\prime}-B\right|_{j\left(l\right)} & \le\sum_{j\in \cJ}\left| (B^{\prime}-B)\right| _{j}\sum_{l=1}^{m_{j}}\frac{\kappa_{l}}{\sqrt{l}} \\ &\le \sqrt{\sum_{j=1}^{d} \left(\sum_{l=1}^{m_{j}}\frac{\kappa_{l}}{\sqrt{l}}\right)^2} \sqrt{\sum_{j^{\prime}=1}^{d_0} | B^{\prime} - B |_{(j^{\prime})}^{2}} \\
&\le \frac{1}{\nu_{gS}(d_0)} \sqrt{\sum_{j=1}^{d} \left(\sum_{l=1}^{m_{j}}\frac{\kappa_{l}}{\sqrt{l}}\right)^2} \| V^{1/2}(B^{\prime}-B) \|_F
\end{aligned} \]

Thus, for any \(B - B^{\prime} \in \mathcal{T}(\rho^*)\), we found \(Z \in \partial \| \cdot \|_{\kappa,\lambda}(B)\) such that 
\[
tr\left(Z^\top\left(B^{\prime}-B\right)\right)\ge\rho^*-2\frac{1}{\nu_{gS}(d_0)} \left(\sum_{j=1}^{d_{0}}\frac{\lambda_{j}}{\sqrt{j}}+\sqrt{\sum_{j=1}^{d} \left(\sum_{l=1}^{m_{j}}\frac{\kappa_{l}}{\sqrt{l}}\right)^2 }\right) \| V^{1/2}(B^{\prime}-B) \|_F.
\]
Hence,
\[
\begin{split}
\sparsity(\rho^*) &= \inf_{B^{\prime\prime} \in \mathcal{T}(\rho^*)} \sup_{Z \in \partial \| \cdot \|_{\kappa,\lambda}(B)} tr\left(Z^\top B^{\prime\prime}\right) \\
&\ge\rho^*-2\left(\sum_{j=1}^{d_{0}}\frac{\lambda_{j}}{\sqrt{j}}+\sqrt{\sum_{j=1}^{d} \left(\sum_{l=1}^{m_{j}}\frac{\kappa_{l}}{\sqrt{l}}\right)^2 }\right) \sqrt{\frac{C_0 Rad(\cB_{\lambda, \kappa}) C}{\nu_{gS}(d_0) \sqrt{n}} \rho^*}.
\end{split}
\]
and, therefore, for $\rho^*$ from (\ref{eq:rho_star_sgs}),  \(\sparsity(\rho^*) \ge \frac{4}{5}\rho^*\).

\subsection{Proof of Lemma \ref{lem:rademacher_sgs}}

To prove Lemma \ref{lem:rademacher_sgs} we first bound the empirical Rademacher complexity \(\widehat{Rad}(\cB_{\kappa, \lambda})\). As a first step, we bound the empirical Rademacher complexity by the empirical Gaussian complexity
$$
\widehat{G}(\cB_{\kappa,\lambda})=\E_G\left\{\frac{1}{\sqrt n} \sup_{B \in \cB_{\kappa,\lambda}}\sum_{i=1}^n\sum_{l=1}^L G_{il}\bbeta_l^T  \bX_i\Big|\bX_1=\bx_1,\ldots,\bX_n=\bx_n\right\}=\E_G\left\{\frac{1}{\sqrt n} \sup_{B \in \cB_{\kappa,\lambda}} tr(B^{\top} Z)\right\},
$$
where $G_{il}$ are i.i.d. $N(0,1)$ and $Z=X^T G$.
We have $\widehat{Rad}(\cB_{\kappa, \lambda}) \leq \sqrt{\frac{\pi}{2}} \widehat{G}(\cB_{\kappa,\lambda})$ \citep[see, e.g., ][Section 5.2]{wainwright2019highdimensional}.

\ignore{
Let \(G \in \R^{n \times L}\) such that \(G_{jl} \sim \mathcal{N}(0,1)\) and denote \(Z = \frac{1}{\sqrt{n}} X^{\top} G\). Then, 
\[
\begin{split}
Rad(\cB_{\kappa, \lambda}) &= \E_{X}\left\{\widehat{Rad}(\cB_{\kappa, \lambda})\right\} = \E_{X} \E_{\Sigma}\left\{ \sup_{B\in \cB_{\kappa, \lambda}} tr\left(\frac{1}{\sqrt{n}} \Sigma^{\top} X B\right)\right\} \\
&\le \sqrt{\frac{\pi}{2}}\E_{X} \E_{G}\left\{ \sup_{B\in \cB_{\kappa, \lambda}} tr\left(Z^{\top} B\right)\right\}.
\end{split}
\]
}

Define 
$$
\delta_j=\sqrt{\sum_{l=1}^{L}\left(\left|Z_{j(l)}\right|- \sqrt{\frac{2}{\pi}}\frac{7}{1440} |X|_{2j}~ \kappa_{l}\right)_{+}^{2}},\;\;j=1,\ldots,d.
$$
To bound $\widehat{G}(\cB_{\kappa,\lambda})$ we need the following two lemmas: 

\begin{lemma} \label{lem:rad1}
\begin{align*}
\widehat{G}(\cB_{\kappa,\lambda}) \leq		\sqrt{\frac{2}{\pi}}\frac{7}{1440} \max_{1 \leq j \leq d} |X_{\cdot j}|_2 +\E_G \max_{1 \leq j \leq d}\frac{\delta_{(j)}}{\lambda_j}.
\end{align*}
\end{lemma}

\begin{lemma} \label{lem:rad_truncated} Let \(\kappa_{L} \ge \sqrt{\frac{\pi}{2}}\frac{2880}{7 \sqrt{n}}\). Then,
	conditionally on $X$,
	$$
	\E_{G} \max_{1 \leq j \leq d}\frac{\delta_{(j)}}{\lambda_j} \leq  C_{0} \max_{1\le j\le d} \left\{\frac{1}{\sqrt{n}} | X|_{2j}~ \frac{\sqrt{2 \sum_{j=1}^{L}\frac{1}{l}\left(\frac{Le}{l}\right)^{l}e^{-C^2 n l \kappa_{l}^{2}}+2\log\left(\frac{de}{j}\right)}}{\lambda_{j}}\right\},
    $$
	where \(C = \sqrt{\frac{2}{\pi}}\frac{7}{2880}\) and \(C_0 = 2\left(1 + \sqrt{\pi}\right)\).
\end{lemma}

Lemmas \ref{lem:rad1} and \ref{lem:rad_truncated} together imply
\[
\widehat{Rad}(\cB_{\kappa, \lambda}) \le \left(\frac{7}{1440 }\sqrt{n} + C_{0} \sqrt{\frac{\pi}{2}}\max_{1\le j \le d} \frac{\sqrt{2 \sum_{j=1}^{L}\frac{1}{l}\left(\frac{Le}{l}\right)^{l}e^{-C^2 n l \kappa_{l}^{2}}+2\log\left(\frac{de}{j}\right)}}{\lambda_{j}}\right) \max_{1 \leq j \leq d} \frac {1}{\sqrt{n}} |X_{\cdot j}|_2.
\]
Hence, 
\[
\begin{split}
Rad(\cB_{\kappa, \lambda}) &= \E_{X}\left\{\widehat{Rad}(\cB_{\kappa, \lambda})\right\} \\
&\le \frac{7}{1440} \sqrt{n} +  C_{0} \sqrt{\frac{\pi}{2}}\max_{1\le j \le d} \frac{\sqrt{8 \sum_{j=1}^{L}\frac{1}{l}\left(\frac{Le}{l}\right)^{l}e^{-C^2 n l \kappa_{l}^{2}}+2\log\left(\frac{de}{j}\right)}}{\lambda_{j}}~.
\end{split}
\]

\subsubsection*{Proof of Lemma \ref{lem:rad1}}
	Define two unit balls w.r.t. $||\cdot||_{\kappa,0}$ and $||\cdot||_{0,\lambda}$: \(\mathcal{B}_\kappa=\left\{ B:\sum_{j=1}^d \sum_{l=1}^{L}\kappa_l\left|B_{j(l)}\right|\le1\right\}\)
	and \(\mathcal{B}_{\lambda}=\left\{ B:\sum_{j=1}^{d}\lambda_j|B|_{(j)} \le1\right\} \) and note
	that \(\mathcal{B}_{\kappa,\lambda}\subseteq \mathcal{B}_{\kappa}\cap\mathcal{B}_{\lambda}\). 
	
	For any matrix $A \in\mathbb{R}^{d\times L}$ we have
	\begin{align*}
	E_{G}\sup_{B\in\mathcal{B}_{\kappa,\lambda}}\langle Z,B\rangle
	& \le \E_{G}\sup_{B\in\mathcal{B}_{\kappa}\cap\mathcal{B}_{\lambda}}\langle Z,B\rangle 
		~=~\E_{G} \sup_{B\in\mathcal{B}_{\kappa}\cap\mathcal{B}_{\lambda}}\left\{\left\langle A,B\right\rangle +\left\langle Z-A,B\right\rangle\right\}\\
		& \le \E_{G}\left\{\sup_{B\in\mathcal{B}_{\kappa}\cap\mathcal{B}_{\lambda}}\left\langle A,B\right\rangle +\sup_{B\in\mathcal{B}_{\kappa}\cap\mathcal{B}_{\lambda}}\left\langle Z-A,B\right\rangle \right\}\\
		& \le \E_{G} \left\{\sup_{B\in\mathcal{B}_{\kappa}}\left\langle A,B\right\rangle +\sup_{B\in\mathcal{B}_{\lambda}}\left\langle Z-A,B\right\rangle \right\}.
	\end{align*}
	Similar to the results  for the group Slope of \citet{abramovich2021multiclass},  \(\sup_{B\in\mathcal{B}_{\kappa}}\left\langle A,B\right\rangle  \le \max_{jl}\frac{\left|A_{j\left(l\right)}\right|}{\kappa_{l}}\) and \(\sup_{B\in\mathcal{B}_{\lambda}}\left\langle Z-A,B\right\rangle  \le \max_j\frac{|Z-A|_{(j)}}{\lambda_j}\). 
	Thus, 
	\begin{align*}
	    E\sup_{B\in\mathcal{B}_{\kappa,\lambda}}\langle Z,B\rangle \le \E\left\{\max_{j,l}\frac{\left|A_{j\left(l\right)}\right|}{\kappa_{l}}+\max_{j}\frac{|Z-A|_{(j)}}{\lambda_j}\right\}.
	\end{align*}
	In particular, consider a matrix $A$ such that \(A_{j\left(l\right)}=\text{sign}\left(Z_{j\left(l\right)}\right)\min\left\{ \left|Z_{j\left(l\right)}\right|, \sqrt{\frac{2}{\pi}} \frac{7}{1440} |X_{\cdot j}|_2 \kappa_{l}\right\} \). We then have
	\begin{align*}
	     E\sup_{B\in\mathcal{B}_{\kappa,\lambda}}\langle Z,B\rangle & \le \sqrt{\frac{2}{\pi}}\frac{7}{1440} \max_{1 \leq j \leq d} |X_{\cdot j}|_2 +\E\left\{\max_{1 \leq j \leq d}\frac{\left(\sqrt{\sum_{l=1}^{L}\left(\left|Z_{j(l)}\right|- \sqrt{\frac{2}{\pi}}\frac{7}{1440} |X_{\cdot j}|_2~ \kappa_{l}\right)_{+}^{2}}\right)_{\left(j\right)}}{\lambda_{j}}\right\}.
	\end{align*}

\subsubsection*{Proof of Lemma \ref{lem:rad_truncated}}
	Denoting \(C = \frac{1}{\sqrt{2\pi}} \frac{7}{1440}\), we have
	\begin{equation} \label{eq:work3}
	\begin{split}
		\E_{G}\left\{\left(\frac{\left|Z_{j\left(l\right)}\right|}{\frac{1}{\sqrt{n}}|X_{\cdot j}|_2}- 2C\sqrt{n} \kappa_{l}\right)_{+}^{2}\right\} &=\int_{0}^{\infty}2sP\left(\left(\frac{\left|Z_{j\left(l\right)}\right|}{ \frac{1}{\sqrt{n}}| X_{\cdot j}|_2}- 2C \sqrt{n} \kappa_{l}\right)_{+}^{2}>s^{2}\right)ds\\
		& \le\int_{0}^{\infty}2sP\left(\frac{\left|Z_{j\left(l\right)}\right|}{ \frac{1}{\sqrt{n}} | X_{\cdot j}|_2}>s+ 2C \sqrt{n}\kappa_{l}\right)ds\\
		& \le\int_{0}^{\infty}2s{\binom{L}{l}}P\left(\frac{\left|Z_{jl}\right|}{ \frac{1}{\sqrt{n}} | X_{\cdot j}|_2}>s+ 2C \sqrt{n} \kappa_{l}\right)^{l}ds.
	\end{split}
	\end{equation}
	Note that conditionally on $X$, \(\frac{Z_{jl}}{ \frac{1}{\sqrt{n}} |X_{\cdot j}|_2} \) is an \(\mathcal{N}(0,1)\) Gaussian random variable and, therefore,
	(\ref{eq:work3}) yields
	\be \label{eq:work4}
	\begin{split}
		 \E_{G}\left\{\left(\frac{\left|Z_{j\left(l\right)}\right|}{ \frac{1}{\sqrt{n} } |X_{\cdot j}|_2}- 2C \sqrt{n}  \kappa_{l}\right)_{+}^{2}\right\} & \le\int_{0}^{\infty} 2^{l+1} s{\binom{L}{l}}e^{-\frac{l\left(s+2C\sqrt{n} \kappa_{l}\right)^{2}}{2}}ds \\
		& \le\int_{0}^{\infty}2^{l+1}\frac{1}{l}{\binom{L}{l}}l\left(s+2C\sqrt{n} \kappa_{l}\right)e^{-\frac{l\left(s+ 2C\sqrt{n} \kappa_{l}\right)^{2}}{2}}ds\\
		& =2^{l+1}\frac{1}{l}{\binom{L}{l}}e^{-2C^2 nl\kappa_{l}^{2}}\le 2^{l+1}\frac{1}{l}\left(\frac{Le}{l}\right)^{l}e^{-2C^2 nl\kappa_{l}^{2}}.
	\end{split}
	\ee
	For \( \kappa_{l} \ge \frac{1}{C\sqrt{n}}\), (\ref{eq:work4}) implies
	\[
	\E \left\{\left(\left|Z_{j\left(l\right)}\right|- 2C |X_{\cdot j}|_2 \kappa_{l}\right)_{+}^{2}\right\} \le 2\frac{1}{l}\left(\frac{Le}{l}\right)^{l}e^{-C^2 nl\kappa_{l}^{2}}~ \frac{1}{n} | X_{\cdot j}|_2^2
	\]
	Hence, by Jensen inequality,
	\[
	\begin{split}
	\mathbb{E}\sqrt{\sum_{l=1}^{L}\left(\left|Z_{j\left(l\right)}\right|- 2C | X_{\cdot j}|_2 \kappa_l\right)_{+}^{2}} &\le \sqrt{\mathbb{E}\left[\sum_{l=1}^{L}\left(\left|Z_{j\left(l\right)}\right|- 2C | X_{\cdot j}|_2 \kappa_{l}\right)_{+}^{2}\right]} \\
	&\le \sqrt{2 \sum_{j=1}^{L}\frac{1}{l}\left(\frac{Le}{l}\right)^{l}e^{-C^2 n l \kappa_{l}^{2}}}~ \frac{1}{\sqrt{n}}|X_{\cdot j}|_2.
	\end{split}
	\]
	Let 
	\[
	M_{j} = \sqrt{2 \sum_{j=1}^{L}\frac{1}{l}\left(\frac{Le}{l}\right)^{l}e^{-C^2 n l \kappa_{l}^{2}}}.
	\]
	
	One can verify that the function \(f_j({\bf z}) = \sqrt{\sum_{l=1}^{L} \left(|z_{j(l)}| - 2C |X_{\cdot j}|_{2} \kappa_{l}\right)_{+}^2}: \mathbb{R}^L \rightarrow \mathbb{R}\) is a 1-Lipschitz function. Recall that \(Z_{jl} \sim \mathcal{N}(0,\frac{1}{\sqrt{n}}|X_{\cdot j}|_2)\) and,
	therefore, by the Tsirelson-Ibragimov-Sudakov inequality \citep[Theorem 5.6]{boucheron2013concentration}, for any \(s,u \ge 1\),
	\[
	\begin{split}
	P&\left( f_{j}(Z) > s \frac{1}{\sqrt{n}}|X_{\cdot j}|_2 \sqrt{2M_{j}^2 + 2u} \right) \le 
	P\left( f_{j}(Z) >  \frac{1}{\sqrt{n}}|X_{\cdot j}|_2 M_{j} + \frac{1}{\sqrt{n}}|X_{\cdot j}|_2 s\sqrt{u}\right)  \\ 
	\le& P\left( f_{j}(Z) >  E f_{j}(Z) + \frac{1}{\sqrt{n}}|X_{\cdot j}|_2 s\sqrt{u} \right) \le e^{-\frac{s^2}{2}u}.
	\end{split}
	\]
	
	Thus, for $s\ge 2,$ we have,
	\[
	\begin{split}
	P&\left(\frac{f_{\left(j\right)}(Z)}{\lambda_{j}} >s \frac{1}{\sqrt{n}} |X_{\cdot j}|_2 \frac{\sqrt{2 M^{2}_{j}+2 \log\left(de/j\right)}}{\lambda_{j}}\right) \\
	&\le \binom{d}{j} P\left(\frac{f_{j}(Z)}{\lambda_{j}} >s \frac{1}{\sqrt{n}} |X_{\cdot j}|_2 \frac{\sqrt{2 M^{2}_{j}+2 \log\left(de/j\right)}}{\lambda_{j}}\right)^{j} \\
	& \le \binom{d}{j} e^{-j \frac{s^2}{2} \log\left(de/j\right)} \le \left(\frac{de}{j}\right)^{-j\left(\frac{s^{2}}{2}-1\right)}\le \left(\frac{de}{j}\right)^{-j\frac{s^{2}}{4}},
	\end{split}
	\]
	and applying the union bound, 
	\begin{equation} \label{eq:work5}
		\begin{split}
			P\left(\max_{j} \frac{f_{\left(j\right)}(Z)}{\lambda_{j}}>s\max_{j} \frac{1}{\sqrt{n}} |X_{\cdot j}|_2 \frac{\sqrt{2 M^{2}_{j}+2\log\left(de/j\right)}}{\lambda_{j}}\right) &\le      \sum_{j=1}^{d}\left(\frac{de}{j}\right)^{-j\frac{s^{2}}{4}} \le \sum_{j=1}^{d}e^{-j\frac{s^{2}}{4}} \\
			&\le \frac{e^{-\frac{s^{2}}{4}}}{1-e^{-\frac{s^{2}}{4}}} \le 2 e^{-\frac{s^{2}}{4}},
		\end{split}
	\end{equation}
	Finally, (\ref{eq:work5}) implies
	\begin{equation*}
		\begin{split}
			\mathbb{E}&\left\{\frac{\max_{1 \leq j \leq d}\frac{\delta_{(j)}}{\lambda_j}}
			{\max_{j=1}^{d} \frac{1}{\sqrt{n}} | X_{\cdot j}|_{2} \frac{\sqrt{2 \sum_{j=1}^{L}\frac{1}{l}\left(\frac{Le}{l}\right)^{l}e^{-C^2 n l \kappa_{l}^{2}}+2\log\left(\frac{de}{j}\right)}}{\lambda_{j}}}\right\} \\
			&= \int_{0}^{\infty}P\left(\max_{j}\frac{f_{\left(j\right)}(Z)}{\lambda_{j}}>s\max_{j} \frac{1}{\sqrt{n}} |X_{\cdot j}|_2 \frac{\sqrt{2 M^{2}_{j}+2\log\left(de/j\right)}}{\lambda_{j}}\right) \\ &\le 2\left(1 + \sqrt{\pi}\right).
		\end{split}
	\end{equation*}

\subsection{Proof of Lemma \ref{lem:nu_rademacher}}

Let \(U \in \R^{L \times (L-1)}\) be a matrix with orthonormal columns such that \(U U^\top = I - \frac{1}{L} \bone \bone^{\top}\). One can easily verify that $B=BUU^\top$. Recall that
\[
Rad(\cB_{\lambda}) = \E_{X} \E_{\Sigma} \left[\frac{1}{\sqrt{n}} \sup_{B \in \cB_{\lambda}} tr(\Sigma U U^{\top} B^{\top} X^{\top})\right] =  \E_{X} \E_{\Sigma} \left[\frac{1}{\sqrt{n}} \sup_{B \in \cB_{\lambda}} tr(U^{\top} B^{\top} K)\right],
\]
where \(K = X^{\top} \Sigma U \in \R^{d\times (L - 1)}\). By duality of Schatten norms,
\[
\frac{1}{\sqrt{n}} \sup_{B \in \cB_{\lambda}} tr(U^{\top} B^{\top} K) =
\frac{1}{\lambda} \frac{1}{\sqrt{n}} \sup_{\|B\|_* \leq 1} tr(U^{\top} B^{\top} K)
=\frac{1}{\lambda}\frac{1}{\sqrt{n}} \| X^{\top} \Sigma U \|_{2}.
\]

Denote \(v(X) = \|\frac{1}{\sqrt{n}} X\|_{2}\) and 
\(\omega(X)=\|\frac{1}{\sqrt{n}} X\|_F \leq v(X) \sqrt{d}\). 
\ignore{
and note that 
\[
\left\| \frac{1}{\sqrt{n}} K \right\|_{F} = \left\| \frac{1}{\sqrt{n}} X^T \Sigma \right\|_{F} = \sqrt{\frac{1}{n} \sum_{i=1}^{n} \sum_{j=1}^{d} X_{ij}^2} = 1.
\]
}
By Theorem 3.2 of \citet{rudelson2013hansonwright}, conditionally on $X$, for any $s, t>1$   
\begin{equation} \label{eq:rud-ver}
P\left(\frac{1}{\sqrt{n}} \| X^{\top} \Sigma U \|_{2} > C\left(s \omega(X)+t\sqrt{L-1}~ v(X)\right) \Big| X\right) \le 2 \exp\left(-\frac{\omega^2(X)}{v^2(X)}s^2 - (L-1) t^2\right),
\end{equation}
where $C>0$ is given in their theorem.

Assume first that  $v(X)\sqrt{L-1} \geq \omega(X)$. Take \(s = \frac{v(X)}{\omega(X)} \sqrt{L-1} t>1\) in (\ref{eq:rud-ver}) to get 
\[
P\left(\frac{1}{\sqrt{n}} \| X^{\top} \Sigma \|_{2} > 2C t\sqrt{L-1}~v(X)\right) \le 2 \exp\left(-2t^2 (L-1)\right).
\]
Setting \(u = 2Ct\sqrt{L-1}~ v(X)\) yields 
\[
P\left(\frac{1}{\sqrt{n}} \| X^{\top} \Sigma \|_{2} > u~ \Big|X\right) \le 2 \exp\left(-  \frac{u^2}{2C^2v(X)^2}\right),
\]
for any $u > 2C\sqrt{L-1}v(X)$ and, therefore, the empirical Rademacher complexity 
$$
\widehat{Rad}(\cb_\lambda) \leq \frac{1}{\lambda} \left(2C \sqrt{L-1}v(X)+2\int_{2C \sqrt{L-1}v(X)}^\infty e^{- \frac{u^2}{2C^2v(X)^2}} du \right)  \leq 
C \frac{1}{\lambda} \sqrt{L-1} v(X).
$$

Similarly, for $v(X)\sqrt{L-1} < \omega(X)$, take $t=\frac{\omega(X)}{v(X)\sqrt{L-1}}s>1$ in
(\ref{eq:rud-ver}) and $u=2Cs\omega(X)$ to get
\[
P\left(\frac{1}{\sqrt{n}} \| X^{\top} \Sigma \|_{2} > u~ \Big|X\right) \le 2 \exp\left(-  \frac{u^2}{2C^2v(X)^2}\right),
\]
for any $u>2C\omega(X)$ and, therefore,
$$
\widehat{Rad}(\cb_\lambda) \leq 
C \frac{1}{\lambda}  \omega(X) \leq C \frac{1}{\lambda} v(X) \sqrt{d}.
$$
Combining both cases we have
\begin{equation} \label{eq:rad_empir}
\widehat{Rad}(\cb_\lambda) \leq C \frac{1}{\lambda}v(X)(\sqrt{L-1}+\sqrt{d}).
\end{equation}
To complete the proof of the lemma apply the results of  \citet[Section 5.4.1]{vershynin2012introduction} for sub-Gaussian matrices with independent rows to get
\begin{equation} \label{eq:Ev}
E_Xv(X) \leq \sqrt{E_Xv^2(X)} \leq C \sqrt{\tau_{1}(V)}.
\end{equation}

\section{Proof of Theorem \ref{th:nuclear_lower}} \label{sec:appendix_lower}
Consider the class $\tilde{\cC}_L(r_0)$ of $r_0$-globally sparse linear $L$-class classifiers from Section \ref{subsec:global} but with the {\em known} subset of $r_0$ significant features. Evidently, $\tilde{\cC}_L(r_0) \subset \cC^*_L(r_0)$.  Apply now
the results of \citet[Theorem 2]{abramovich2021multiclass} on the lower bounds for global row-wise sparse classification to get  
\begin{equation} \label{eq:lowerbound1}
\inf_{\widetilde \eta}\sup_{\eta^* \in \cC^*_L(r_0),~\mathbb{P}_X}
\cE(\tilde{\eta}, \eta^{*})  \geq \inf_{\widetilde \eta}\sup_{\eta^* \in \tilde{\cC}_L(r_0),~\mathbb{P}_X}
\cE(\tilde{\eta}, \eta^{*}) \geq C \sqrt{\frac{r_0 (L-1)}{n}}.
\end{equation}

On the other hand, consider $r_0$-class classification, where all $d$ features are significant
($d_0=d$). It is obvious that 
$\cC_{r_0}(d) \subset \cC^*_{r_0}(r_0)$ and
that $r_0$-class classification cannot be harder than the $L$-class one.
Thus, exploiting again Theorem 2 of \citet{abramovich2021multiclass} we have
\begin{equation} \label{eq:lowerbound2}
\inf_{\widetilde \eta}\sup_{\eta^* \in \cC^*_L(r_0),~\mathbb{P}_X}
\cE(\tilde{\eta}, \eta^{*}) \geq
\inf_{\widetilde \eta}\sup_{\eta^* \in \cC^*_{r_0}(r_0),~\mathbb{P}_X}
\cE(\tilde{\eta}, \eta^{*}) \geq
\inf_{\widetilde \eta}\sup_{\eta^* \in \cC_{r_0}(d),~\mathbb{P}_X}
\cE(\tilde{\eta}, \eta^{*})  \geq C \sqrt{\frac{r_0 d}{n}}.
\end{equation}
Combining (\ref{eq:lowerbound1}) and (\ref{eq:lowerbound2}) completes the proof of the theorem.

\section{Sparse group Slope algorithm} \label{sec:sgS_algorithm}

The penalized MLE minimization problem in \eqref{eq:sgS} involves a sum of a convex smooth log-likelihood and a convex but non-smooth penalty consisting of two terms. A common approach  to solve such optimization problems is by the proximal gradient method \citep[e.g.,][]{beck2017first}. A general proximal operator of a given convex function \(f\) is defined as 
\[
\prox_{f}(a) = \arg\min_{b} 
\left\{\frac{1}{2} \|a - b\|^{2} + f(b) \right\}.
\]
For the setup at hand consider the  proximal operator
\be \label{eq:prox}
\prox_{ \| \cdot \|_{\kappa,\lambda} }(A) = \arg\min_{B} 
\left\{\frac{1}{2} \|A - B\|^{2}_{F} + \| B \|_{\kappa,\lambda}\right\},
\ee 
where recall that $ \|B\|_{\kappa,\lambda} =\sum_{j=1}^d \lambda_j|B|_{(j)}+\sum_{j=1}^d\sum_{l=1}^L \kappa_l |B|_{j(l)}=\|B\|_\lambda+\sum_{j=1}^d \|B_{j\cdot}\|_\kappa$.

There exist the efficient proximal gradient descent algorithms for computing proximal operators \(\prox_{\| \cdot \|_{\kappa}}\) and \(\prox_{\| \cdot \|_{\lambda}}\) for \(\| \cdot \|_{\kappa}\) and \(\| \cdot \|_{\lambda}\) separately \citep[see respectively][]{bogdan2015slope, brzyski2019group}. We now show that applying \(\prox_{\| \cdot \|_{\kappa}}\) and \(\prox_{\| \cdot \|_{\lambda}}\) consecutively results in \(\prox_{\| \cdot \|_{
\kappa,\lambda}}\) as depicted by Algorithm \ref{alg:prox}:

\begin{algorithm}[H] \label{alg:prox}
\SetAlgoLined
  \DontPrintSemicolon
  \For{$j \to 1 \ldots d$}{
    \(U_{j \cdot} = \prox_{\| \cdot \|_{\kappa}}(A_{j \cdot}) \)\;
  }
  \(B \leftarrow \prox_{\| \cdot \|_{\lambda}}( U )\) \; 
 \caption{\(\prox_{\| \cdot \|_{\kappa,\lambda}}(A)\)}
\end{algorithm}
The proof relies on the second prox theorem \citep[Theorem 6.39]{beck2017first} and the following general lemma:
\begin{lemma}\label{lem:prox_relation}
Assume that for all \(a\), \(\partial g(\prox_{f}(a)) \supseteq \partial g(a)\), then for all \(b\), \(\prox_{f+g}(b) = \prox_{f}(\prox_g(b))\).
\end{lemma}
\begin{proof}
For a given $b$, let \(a = \prox_{g}(b)\) and \(z = \prox_{f}(a)\). By the second prox theorem, 
\(b - a \in \partial g(a)\) and \(a- z \in \partial f(z)\).
By the condition, \(\partial g(z) \supseteq \partial g(a)\), and therefore,
\[
b - z = b - a + a - z \in \partial f(z) + \partial g(z) = \partial (f + g)(z)
\]
which implies by the second prox theorem that
\(
z = prox_{f+g}(b).
\)
\end{proof}

Applying Lemma \ref{lem:prox_relation} for \(g(A) = \sum_{j=1}^{d} \| A_{{j}, \cdot} \|_{\kappa} \) and \(f(A) = \| A \|_{\lambda} \) relies on the following lemma:
\begin{lemma}\label{lem:gs_subgradient}
For \(Z,A\in \R^{d\times L}\) such that \(Z \in \partial\| \cdot \|_{\lambda}(A)\) and for any \(j \in \{1,\ldots,d\}\), there exists \(c_j\ge0\) such that \(Z_{j \cdot} = c_j A_{j \cdot}\).
\end{lemma}

\begin{proof}
Let \(Z \in \partial\| \cdot \|_{\lambda}(A)\). Thus,
\[
Z \in \argmax_{\| Z \|_{\lambda}^{*} \le 1} tr(Z^{\top}A) = \argmax_{\| Z \|_{\lambda}^{*} \le 1} \sum_{j=1}^{d} Z_{j\cdot}^{\top} A_{j\cdot}~,
\] 
where \(\| \cdot \|_{\lambda}^{*}\) is the dual norm.
Since the norm \(\| \cdot \|_{\lambda}\) is invariant to rotation of the rows, so does its dual norm \(\| \cdot \|_{\lambda}^{*}\) because we can always rotate the rows of the norming matrix. Thus, the maximum above is when \(Z_{j\cdot} = c_{j}A_{j\cdot}\) for some \(c_j \ge 0\).
\end{proof}

Let \(Z = \prox_{\| \cdot \|_{\lambda}}(A)\). By the second prox theorem  we have \(A - Z \in \partial \| \cdot \|_{\lambda}(Z)\), and by Lemma \ref{lem:gs_subgradient}, \(A_{j} - Z_{j} = c_{j} Z_{j}\) for some \(c_j > 0\).  Thus, \(Z_{j} = \frac{1}{1 + c_j} A_j\). 

Let $V \in \partial (\sum_{j=1}^{d} \| e_j^{\top} \cdot\|_{\kappa})(A)$, that is, \(V_{j \cdot} \in \partial \| \cdot \|_{\kappa} (A_{j\cdot})\). By the definition of the subgradient, for any $\bu \in \R^{L}$,
\[
    \| A_{j \cdot} \|_{\kappa} + V^{\top}_{j \cdot}(\bu - A_{j\cdot}) \le \| \bu \|_{\kappa}  
\]
Let \(\bu^{\prime} \in \R^{L}\). Then,
\[
\begin{split}
    \| Z_{j \cdot} \|_{\kappa} + V^{\top}_{j \cdot}(\bu^{\prime} - Z_{j\cdot}) &= \frac{1}{1 + c_j}\| A_{j \cdot} \|_{\kappa} + \frac{1}{1 + c_j} V^{\top}_{j \cdot}((1+c_j)\bu^{\prime} - A_{j\cdot}) \\
    &\le \frac{1}{1 + c_j} \| (1+c_j)\bu^{\prime} \|_{\kappa} = \| \bu^{\prime} \|_{\kappa}
\end{split}
\]
and, therefore, \(
V^{\top}_{j \cdot}(\bu^{\prime} - Z_{j\cdot}) \le \| \bu^{\prime} \|_{\kappa} - \| Z_{j \cdot} \|_{\kappa}\) implying \(V_{j\cdot} \in \partial \| \cdot \|_{\kappa}(Z_{j\cdot})\). Hence, \(V \in \partial (\sum_{j=1}^{d} \| e_j^{\top} \cdot\|_{\kappa})(Z)\) and the condition for Lemma \ref{lem:prox_relation} holds, i.e. 
\[
\partial \| \cdot \|_{\kappa}(\prox_{\|\cdot \|_{\lambda}}(A)) \supseteq \partial \| \cdot \|_{\kappa}(A).
\]

\newpage
\bibliography{reference}

\begin{thebibliography}{40}
\providecommand{\natexlab}[1]{#1}
\providecommand{\url}[1]{\texttt{#1}}
\expandafter\ifx\csname urlstyle\endcsname\relax
  \providecommand{\doi}[1]{doi: #1}\else
  \providecommand{\doi}{doi: \begingroup \urlstyle{rm}\Url}\fi

\bibitem[Abramovich and Grinshtein(2019)]{abramovich2019highdimensional}
Felix Abramovich and Vadim Grinshtein.
\newblock High-{{dimensional classification}} by {{sparse logistic
  regression}}.
\newblock \emph{IEEE Transactions on Information Theory}, 65\penalty0
  (5):\penalty0 3068--3079, May 2019.
\newblock ISSN 0018-9448, 1557-9654.
\newblock \doi{10.1109/TIT.2018.2884963}.

\bibitem[Abramovich et~al.(2021)Abramovich, Grinshtein, and
  Levy]{abramovich2021multiclass}
Felix Abramovich, Vadim Grinshtein, and Tomer Levy.
\newblock Multiclass {{classification}} by {{sparse multinomial logistic
  regression}}.
\newblock \emph{IEEE Transactions on Information Theory}, 67\penalty0
  (7):\penalty0 4637--4646, July 2021.
\newblock ISSN 1557-9654.
\newblock \doi{10.1109/TIT.2021.3075137}.

\bibitem[Alquier et~al.(2019)Alquier, Cottet, and
  Lecu{\'e}]{alquier2019estimation}
Pierre Alquier, Vincent Cottet, and Guillaume Lecu{\'e}.
\newblock Estimation bounds and sharp oracle inequalities of regularized
  procedures with {{Lipschitz}} loss functions.
\newblock \emph{The Annals of Statistics}, 47\penalty0 (4):\penalty0
  2117--2144, August 2019.
\newblock ISSN 0090-5364, 2168-8966.
\newblock \doi{10.1214/18-AOS1742}.

\bibitem[Bach(2008)]{bach2008consistency}
Francis~R. Bach.
\newblock Consistency of {{trace norm minimization}}.
\newblock \emph{The Journal of Machine Learning Research}, 9:\penalty0
  1019--1048, June 2008.
\newblock ISSN 1532-4435.

\bibitem[Beck(2017)]{beck2017first}
Amir Beck.
\newblock \emph{First-order Methods in Optimization}.
\newblock SIAM, 2017.

\bibitem[Bellec et~al.(2018)Bellec, Lecu{\'e}, and Tsybakov]{bellec2018slope}
Pierre~C. Bellec, Guillaume Lecu{\'e}, and Alexander~B. Tsybakov.
\newblock Slope meets {{Lasso}}: {{improved}} oracle bounds and optimality.
\newblock \emph{The Annals of Statistics}, 46\penalty0 (6B):\penalty0
  3603--3642, December 2018.
\newblock ISSN 0090-5364, 2168-8966.
\newblock \doi{10.1214/17-AOS1670}.

\bibitem[Bickel and Levina(2004)]{bickel2004fisher}
Peter~J. Bickel and Elizaveta Levina.
\newblock {Some theory for Fisher's linear discriminant function, `naive
  Bayes', and some alternatives when there are many more variables than
  observations}.
\newblock \emph{Bernoulli}, 10\penalty0 (6):\penalty0 989 -- 1010, 2004.
\newblock \doi{10.3150/bj/1106314847}.
\newblock URL \url{https://doi.org/10.3150/bj/1106314847}.

\bibitem[Bickel et~al.(2009)Bickel, Ritov, and
  Tsybakov]{bickel2009simultaneous}
Peter~J. Bickel, Ya'acov Ritov, and Alexander~B. Tsybakov.
\newblock Simultaneous analysis of {{Lasso}} and {{Dantzig}} selector.
\newblock \emph{The Annals of Statistics}, 37\penalty0 (4):\penalty0
  1705--1732, August 2009.
\newblock ISSN 0090-5364.
\newblock \doi{10.1214/08-AOS620}.

\bibitem[Bogdan et~al.(2015)Bogdan, Van Den~Berg, Sabatti, Su, and
  Cand{\`e}s]{bogdan2015slope}
Ma{\l}gorzata Bogdan, Ewout Van Den~Berg, Chiara Sabatti, Weijie Su, and
  Emmanuel~J Cand{\`e}s.
\newblock Slope—adaptive variable selection via convex optimization.
\newblock \emph{The Annals of Applied Statistics}, 9\penalty0 (3):\penalty0
  1103, 2015.

\bibitem[Boucheron et~al.(2013)Boucheron, Lugosi, and
  Massart]{boucheron2013concentration}
St{\'e}phane Boucheron, G{\'a}bor Lugosi, and Pascal Massart.
\newblock \emph{Concentration {{Inequalities}}: {{A Nonasymptotic Theory}} of
  {{Independence}}}.
\newblock {OUP Oxford}, February 2013.
\newblock ISBN 978-0-19-953525-5.

\bibitem[Brzyski et~al.(2019)Brzyski, Gossmann, Su, and
  Bogdan]{brzyski2019group}
Damian Brzyski, Alexej Gossmann, Weijie Su, and Ma{\l}gorzata Bogdan.
\newblock Group slope--adaptive selection of groups of predictors.
\newblock \emph{Journal of the American Statistical Association}, 114\penalty0
  (525):\penalty0 419--433, 2019.

\bibitem[Bunea et~al.(2011)Bunea, She, and Wegkamp]{bunea2011optimal}
Florentina Bunea, Yiyuan She, and Marten~H. Wegkamp.
\newblock Optimal selection of reduced rank estimators of high-dimensional
  matrices.
\newblock \emph{The Annals of Statistics}, 39\penalty0 (2):\penalty0
  1282--1309, April 2011.
\newblock ISSN 0090-5364, 2168-8966.
\newblock \doi{10.1214/11-AOS876}.

\bibitem[Candes and Plan(2010)]{candes2010matrix}
Emmanuel~J. Candes and Yaniv Plan.
\newblock Matrix {{completion with noise}}.
\newblock \emph{Proceedings of the IEEE}, 98\penalty0 (6):\penalty0 925--936,
  June 2010.
\newblock ISSN 1558-2256.
\newblock \doi{10.1109/JPROC.2009.2035722}.

\bibitem[Chen and Sun(2006)]{chen2006consistency}
Di-Rong Chen and Tao Sun.
\newblock Consistency of {{multiclass empirical risk minimization methods
  based}} on {{convex loss}}.
\newblock \emph{The Journal of Machine Learning Research}, 7\penalty0
  (86):\penalty0 2435--2447, 2006.
\newblock ISSN 1533-7928.

\bibitem[Chen and Lee(2021)]{chen2021binary}
Le-Yu Chen and Sokbae Lee.
\newblock Binary classification with covariate selection through
  {$\ell_0$}-penalised empirical risk minimisation.
\newblock \emph{The Econometrics Journal}, 24\penalty0 (1):\penalty0 103--120,
  January 2021.
\newblock ISSN 1368-4221.
\newblock \doi{10.1093/ectj/utaa017}.

\bibitem[Daniely et~al.(2012)Daniely, Sabato, and
  Shwartz]{daniely2012multiclass}
Amit Daniely, Sivan Sabato, and Shai Shwartz.
\newblock Multiclass {{learning approaches}}: {{a theoretical comparison}} with
  {{implications}}.
\newblock In \emph{Advances in {{Neural Information Processing Systems}}},
  volume~25. {Curran Associates, Inc.}, 2012.

\bibitem[Daniely et~al.(2015)Daniely, Sabato, {Ben-David}, and
  {Shalev-Shwartz}]{daniely2015multiclass}
Amit Daniely, Sivan Sabato, Shai {Ben-David}, and Shai {Shalev-Shwartz}.
\newblock Multiclass {{learnability}} and the {{ERM principle}}.
\newblock \emph{The Journal of Machine Learning Research}, 16\penalty0
  (72):\penalty0 2377--2404, 2015.
\newblock ISSN 1533-7928.

\bibitem[Fan and Fan(2008)]{fan2008high}
Jianqing Fan and Yingying Fan.
\newblock {High-dimensional classification using features annealed independence
  rules}.
\newblock \emph{The Annals of Statistics}, 36\penalty0 (6):\penalty0 2605 --
  2637, 2008.
\newblock \doi{10.1214/07-AOS504}.
\newblock URL \url{https://doi.org/10.1214/07-AOS504}.

\bibitem[Friedman et~al.(2010)Friedman, Hastie, and
  Tibshirani]{friedman2010regularization}
Jerome Friedman, Trevor Hastie, and Rob Tibshirani.
\newblock Regularization {{paths}} for {{generalized linear models}} via
  {{coordinate descent}}.
\newblock \emph{Journal of Statistical Software}, 33\penalty0 (1):\penalty0
  1--22, 2010.
\newblock ISSN 1548-7660.

\bibitem[Koltchinskii and Panchenko(2002)]{koltchinskii2002empirical}
Vladimir Koltchinskii and Dmitry Panchenko.
\newblock Empirical {{margin distributions}} and {{bounding}} the
  {{generalization error}} of {{combined classifiers}}.
\newblock \emph{The Annals of Statistics}, 30\penalty0 (1):\penalty0 1--50,
  2002.
\newblock ISSN 0090-5364.

\bibitem[Koltchinskii et~al.(2011)Koltchinskii, Lounici, and
  Tsybakov]{koltchinskii2011nuclearnorm}
Vladimir Koltchinskii, Karim Lounici, and Alexander~B. Tsybakov.
\newblock Nuclear-norm penalization and optimal rates for noisy low-rank matrix
  completion.
\newblock \emph{The Annals of Statistics}, 39\penalty0 (5):\penalty0
  2302--2329, October 2011.
\newblock ISSN 0090-5364, 2168-8966.
\newblock \doi{10.1214/11-AOS894}.

\bibitem[Lecu{\'e} and Mendelson(2018)]{lecue2018regularization}
Guillaume Lecu{\'e} and Shahar Mendelson.
\newblock Regularization and the small-ball method {{I}}: {{sparse}} recovery.
\newblock \emph{The Annals of Statistics}, 46\penalty0 (2):\penalty0 611--641,
  April 2018.
\newblock ISSN 0090-5364, 2168-8966.
\newblock \doi{10.1214/17-AOS1562}.

\bibitem[Lei et~al.(2019)Lei, Dogan, Zhou, and Kloft]{lei2019datadependent}
Yunwen Lei, {\"U}r{\"u}n Dogan, Ding-Xuan Zhou, and Marius Kloft.
\newblock Data-{{dependent generalization bounds}} for {{multi-class
  classification}}.
\newblock \emph{IEEE Transactions on Information Theory}, 65\penalty0
  (5):\penalty0 2995--3021, May 2019.
\newblock ISSN 1557-9654.
\newblock \doi{10.1109/TIT.2019.2893916}.

\bibitem[Maximov and Reshetova(2016)]{maximov2016tight}
Yurii Maximov and Dar'ya Reshetova.
\newblock Tight risk bounds for multi-class margin classifiers.
\newblock \emph{Pattern Recognition and Image Analysis}, 26\penalty0
  (4):\penalty0 673--680, October 2016.
\newblock ISSN 1555-6212.
\newblock \doi{10.1134/S105466181604009X}.

\bibitem[McCullagh and Nelder(1989)]{mccullagh1989generalized}
Peter McCullagh and John~A Nelder.
\newblock \emph{Generalized Linear Models}.
\newblock Routledge, 1989.

\bibitem[Mohri et~al.(2018)Mohri, Rostamizadeh, and
  Talwalkar]{mohri2018foundations}
Mehryar Mohri, Afshin Rostamizadeh, and Ameet Talwalkar.
\newblock \emph{Foundations of {{Machine Learning}}}.
\newblock Adaptive {{Computation}} and {{Machine Learning}} Series. {MIT
  Press}, {Cambridge, MA, USA}, second edition, December 2018.
\newblock ISBN 978-0-262-03940-6.

\bibitem[Pires and Szepesv{\'a}ri(2016)]{pires2016multiclass}
Bernardo~{\'A}vila Pires and Csaba Szepesv{\'a}ri.
\newblock Multiclass {{classification calibration functions}}.
\newblock \emph{arXiv:1609.06385 [cs, stat]}, September 2016.

\bibitem[Powers et~al.(2018)Powers, Hastie, and Tibshirani]{powers2018nuclear}
Scott Powers, Trevor Hastie, and Robert Tibshirani.
\newblock Nuclear penalized multinomial regression with an application to
  predicting at bat outcomes in baseball.
\newblock \emph{Statistical Modelling}, 18\penalty0 (5-6):\penalty0 388--410,
  2018.

\bibitem[Reeve and Kaban(2020)]{reeve2020optimistic}
Henry Reeve and Ata Kaban.
\newblock Optimistic {{bounds}} for {{multi-output learning}}.
\newblock In \emph{Proceedings of the 37th {{International Conference}} on
  {{Machine Learning}}}, pages 8030--8040. {PMLR}, November 2020.

\bibitem[Rudelson and Vershynin(2013)]{rudelson2013hansonwright}
Mark Rudelson and Roman Vershynin.
\newblock Hanson-{{Wright}} inequality and sub-gaussian concentration.
\newblock \emph{Electronic Communications in Probability}, 18\penalty0
  (none):\penalty0 1--9, January 2013.
\newblock ISSN 1083-589X, 1083-589X.
\newblock \doi{10.1214/ECP.v18-2865}.

\bibitem[{Shalev-Shwartz} and
  {Ben-David}(2014)]{shalev-shwartz2014understanding}
Shai {Shalev-Shwartz} and Shai {Ben-David}.
\newblock \emph{Understanding {{Machine Learning}}: {{From Theory}} to
  {{Algorithms}}}.
\newblock {Cambridge University Press}, {Cambridge}, 2014.
\newblock ISBN 978-1-107-29801-9.
\newblock \doi{10.1017/CBO9781107298019}.

\bibitem[She(2013)]{she2013reduced}
Yiyuan She.
\newblock Reduced rank vector generalized linear models for feature extraction.
\newblock \emph{Statistics and Its Interface}, 6\penalty0 (2):\penalty0
  197--209, 2013.
\newblock ISSN 1938-7997.
\newblock \doi{10.4310/SII.2013.v6.n2.a4}.

\bibitem[Su and Cand{\`e}s(2016)]{su2016slope}
Weijie Su and Emmanuel Cand{\`e}s.
\newblock {{SLOPE}} is adaptive to unknown sparsity and asymptotically minimax.
\newblock \emph{The Annals of Statistics}, 44\penalty0 (3):\penalty0
  1038--1068, June 2016.
\newblock ISSN 0090-5364, 2168-8966.
\newblock \doi{10.1214/15-AOS1397}.

\bibitem[Tsybakov(2004)]{tsybakov2004optimal}
Alexander~B. Tsybakov.
\newblock Optimal aggregation of classifiers in statistical learning.
\newblock \emph{The Annals of Statistics}, 32\penalty0 (1):\penalty0 135--166,
  February 2004.
\newblock ISSN 0090-5364, 2168-8966.
\newblock \doi{10.1214/aos/1079120131}.

\bibitem[van~de Geer(2008)]{geer2008highdimensional}
Sara~A. van~de Geer.
\newblock High-dimensional generalized linear models and the lasso.
\newblock \emph{The Annals of Statistics}, 36\penalty0 (2):\penalty0 614--645,
  April 2008.
\newblock ISSN 0090-5364, 2168-8966.
\newblock \doi{10.1214/009053607000000929}.

\bibitem[Vapnik(2000)]{vapnik2000nature}
Vladimir Vapnik.
\newblock \emph{The Nature of Statistical Learning Theory}.
\newblock Statistics for Engineering and Information Science. Springer, 2000.
\newblock ISBN 978-1-4757-3264-1.

\bibitem[Vershynin(2012)]{vershynin2012introduction}
Roman Vershynin.
\newblock Introduction to the non-asymptotic analysis of random matrices.
\newblock In Gitta Kutyniok and Yonina~C. Eldar, editors, \emph{Compressed
  {{Sensing}}: {{Theory}} and {{Applications}}}, pages 210--268. {Cambridge
  University Press}, {Cambridge}, 2012.
\newblock ISBN 978-1-107-00558-7.
\newblock \doi{10.1017/CBO9780511794308.006}.

\bibitem[Vincent and Hansen(2014)]{vincent2014sparse}
Martin Vincent and Niels~Richard Hansen.
\newblock Sparse group lasso and high dimensional multinomial classification.
\newblock \emph{Computational Statistics $\&$ Data Analysis}, 71\penalty0
  (C):\penalty0 771--786, 2014.
\newblock URL
  \url{https://EconPapers.repec.org/RePEc:eee:csdana:v:71:y:2014:i:c:p:771-786}.

\bibitem[Wainwright(2019)]{wainwright2019highdimensional}
Martin~J. Wainwright.
\newblock \emph{High-{{Dimensional Statistics}}: {{A Non-Asymptotic
  Viewpoint}}}.
\newblock Cambridge {{Series}} in {{Statistical}} and {{Probabilistic
  Mathematics}}. {Cambridge University Press}, {Cambridge}, 2019.
\newblock ISBN 978-1-108-49802-9.
\newblock \doi{10.1017/9781108627771}.

\bibitem[Zhang(2004)]{zhang2004statistical}
Tong Zhang.
\newblock Statistical {{analysis}} of {{some multi-category large margin
  classification methods}}.
\newblock \emph{The Journal of Machine Learning Research}, 5\penalty0
  (Oct):\penalty0 1225--1251, 2004.
\newblock ISSN 1533-7928.

\end{thebibliography}

\end{document}